\documentclass[10pt]{amsart} 
\usepackage{amscd}  
\usepackage{amsfonts}  
\usepackage{amssymb}  
\usepackage{latexsym}
\usepackage{pifont}
 
\newcommand{\ncm}{\newcommand} 
 
\ncm{\beq}{\begin{equation}} 
\ncm{\eeq}{\end{equation}} 
\ncm{\bea}{\begin{eqnarray}} 
\ncm{\eea}{\end{eqnarray}} 
\ncm{\beanon}{\begin{eqnarray*}} 
\ncm{\eeanon}{\end{eqnarray*}} 
 
\newtheorem{thm}{Theorem}[section] 
\newtheorem{pro}[thm]{Proposition} 
\newtheorem{lem}[thm]{Lemma} 
\newtheorem{cor}[thm]{Corollary} 

\theoremstyle{definition}
\newtheorem{lem&def}[thm]{Lemma \& Definition} 
\newtheorem{defi}[thm]{Definition} 
\newtheorem{rmk}[thm]{Remark} 

\theoremstyle{remark}
\newtheorem{exa}[thm]{Example}

\def\Set{\mathsf{Set}} 
 
\def\Ab{\mathsf{Ab}} 
\def\Vec{\mathsf{Vec}}

\def\Mnd{\mathsf{Mnd}}

\def\Cat{\mathsf{Cat}} 
\def\MonCat{\mathsf{MonCat}}
\def\ComonCat{\mathsf{ComonCat}}
\def\Bmd{\mathsf{Bmd}}
\def\Bgd{\mathsf{Bgd}}

\def\LFunc{\mathsf{L}\text{-}\mathsf{Func}}
\def\LMF{\mathsf{L}\text{-}\mathsf{MFunc}}
\def\AMF{\mathsf{A}\text{-}\mathsf{MFunc}} 
\def\ComonFunc{\mathsf{CFunc}}
\def\StMF{\mathsf{st}\text{-}\LMF}
 
\def\Arr{\mathsf{Arr}} 

\def\A{\mathcal{A}}
  \ncm{\oA}{\otimes_\A}  
  \ncm{\iA}{i_\A}        
\def\B{\mathcal{B}} 
  \ncm{\oB}{\otimes_\B}  
  \ncm{\iB}{i_\B}        
\def\C{\mathcal{C}} 
  \ncm{\oC}{\otimes_\C}
  \ncm{\iC}{i_\C}

\def\K{\mathcal{K}} 
\def\L{\mathcal{L}}
\def\M{\mathcal{M}}
  \def\oM{\otimes_\M} \def\oMprime{\otimes_{\M'}} 
  \def\iM{i_\M}       \def\iMprime{i_{\M'}}
\def\N{\mathcal{N}} 
  \def\oN{\otimes_\N} \def\oNprime{\otimes_{\N'}} 
  \def\iN{i_\N}       \def\iNprime{i_{\N'}}
\def\T{\mathcal{T}}

\ncm{\oS}{\,\Box\,}    
\ncm{\iS}{e}
\ncm{\oT}{\otimes}  
\ncm{\iT}{i}

\def\AA{\mathbf{A}}
\def\BB{\mathbf{B}}
\def\CC{\mathbf{C}}
\def\MM{\mathbf{M}}
\def\NN{\mathbf{N}}
\def\OO{\mathbf{O}}

\def\EM{\mathbb{EM}}

\def\twocom{\ast}
\def\onecom{\odot}

\def\bfalpha{\boldsymbol{\alpha}}
\def\bfbeta{\boldsymbol{\beta}}
\def\bfkappa{\boldsymbol{\kappa}}

\ncm{\ttT}{\mathbf{T}}
\ncm{\FT}{F^\mathtt{T}}
\ncm{\ttR}{\mathbf{R}}
\ncm{\ttS}{\mathbf{S}}

\def\ZZ{\mathbb{Z}\,}  

\ncm{\End}{\mbox{\rm End}\,} 
\ncm{\Hom}{\operatorname{Hom}} 
\def\id{\mbox{\rm id}}

\def\to{\rightarrow} 
 
\def\o{\otimes}                   
\def\b{\operatorname{\bullet}}    
\def\c{\operatorname{\circ}}      
\def\x{\times}                 

\ncm{\oalg}{\,\Diamond\,}   
\def\oEM{\oalg}
\def\oTT{\oEM}
\def\ohArr{\,{\scriptstyle\Box}\,}
\def\ovArr{\,{\scriptstyle \blacksquare}\,}
\def\ohBmd{\diamond}

\def\iV{i}  

\newcommand{\dom}{\mathsf{dom}}
\newcommand{\cod}{\mathsf{cod}}
\def\Q{\mathsf{Q}}
\def\QQ{\mathbb{Q}}

\def\bra{\langle} 
\def\ket{\rangle} 
\def\under{\,-\,} 
\ncm{\rarr}[1]{\stackrel{#1}{\longrightarrow}} 
\ncm{\larr}[1]{\stackrel{#1}{\longleftarrow}} 
\def\iso{\rarr{\sim}}
\def\cop{\Delta}

\def\eps{\varepsilon} 
 
\def\du1{\hat 1}

\def\0{_{(0)}} 
\def\1{_{(1)}} 
\def\2{_{(2)}} 
\def\3{_{(3)}}

\def\du1{\hat 1} 
\def\lact{\triangleright} 
 
\def\op{^{\rm op}}
\def\a{\mathbf{a}}
\def\l{\mathbf{l}}
\def\r{\mathbf{r}}
\ncm{\mapsot}{\leftarrow\!\!\!\raisebox{1pt}{$\scriptscriptstyle |$}}

\begin{document} 
 
\title{The monoidal Eilenberg-Moore construction and bialgebroids}
 
\author{Korn\'el Szlach\'anyi}  
\address{Theory Division, Research Institute for Particle and Nuclear Physics,
Budapest, H-1525 Budapest, P. O. Box 49, Hungary}  

\begin{abstract}
Monoidal functors $U\colon \C\to\M$ with left adjoints determine, in a
universal way, monoids $T$ in the category of oplax monoidal endofunctors on
$\M$. Such monads will be called bimonads. Treating bimonads as abstract
"quantum groupoids" we derive Tannaka duality between left adjointable monoidal
functors and bimonads. Bialgebroids, i.e., Takeuchi's
$\times_R$-bialgebras, appear as the special case when $T$ has also a right
adjoint. Street's 2-category of monads then
leads to a natural definition of the 2-category of bialgebroids. 
\end{abstract}

\maketitle

\tableofcontents

\newpage

\section{Introduction}

In the classical theory of monads \cite{MacLane,MacLane-Moerdijk} or triples
 \cite{TTT} one deals with the following construction.  
If a functor $U\colon\C\to\M$ has a left adjoint $F$ with unit $\eta\colon
\M\to UF$ and counit $\eps\colon FU\to\C$ then there is a monad
$\mathtt{T}=\bra UF,U\eps F,\eta\ket$ on $\M$ such that $U$ factorizes
as $U=U^\mathtt{T}K$ through the forgetful functor $U^\mathtt{T}$ of the
category $\M^\mathtt{T}$ of $\mathtt{T}$-algebras with a unique functor
$K\colon \C\to\M^\mathtt{T}$, called the Eilenberg-Moore comparison functor. 
For example, if $A$ is a ring, i.e., a monoid in $\M=\Ab$ then the forgetful
functor $U\colon\,_A\M\to \M$ has a left adjoint, the induction functor
$F=A\o_\ZZ\under$ and the Eilenberg-Moore category $\M^\mathtt{T}$ is
equivalent to the category of left $A$-modules via $K$. 
As it is the case in this example so in general, if $K$ is an equivalence $U$
is called monadic. 

Interpreting the category $\M^\mathtt{T}$ as the category of
$\mathtt{T}$-modules the above construction appears as a primitive
version of Tannaka reconstruction in which adjointable functors are brought
into correspondence with monads on their target categories. Monadicity in
turn plays the role of the representation theorem. Pursuing this idea one can
investigate the extra structure the monad acquires if the forgetful functor is
monoidal. What one obtains is that the monad is oplax monoidal
(called lax comonoidal in this paper) in the sense of the underlying
functor being lax comonoidal and the natural transformations preserve these
comonoidal structures. Such monads will be called \textit{bimonads}.  Bimonads
are the abstract versions of bialgebras in the same spirit as monads are
related to algebras.

Bimonads first appeared in the work of Moerdijk\cite{Mdk} under the name
Hopf monads. He was motivated by generalizing the notion of Hopf operad. In
our context the name bimonad seems to be the more natural as it rhymes with
bialgebra and bialgebroid. Hopf algebroid could then be the name for a bimonad
possessing some sort of antipode.

The motivating example of bimonads is associated to a Takeuchi
$\times_R$-algebra \cite{Takeuchi}, also called \textit{bialgebroid}
\cite{Lu,Brz-Mil,Schauenburg: Duals,Sz}. For a bialgebroid $A$ over $R$ the
algebra $A$ is an extension of $R^e:=R\otimes R\op$, so $A$ is an
$R^e$-$R^e$-bimodule. The underlying functor $T$ of the bimonad is 
$T=A\otimes_{R^e}\under\colon\,_R\M_R\to\,_R\M_R$ where we identified $_R\M_R$ 
with $_{R^e}\M$. Such bimonads are obtained by the above mentioned Tannaka
recontruction from the strict monoidal forgetful functors
$U\colon\,_A\M\to\,_R\M_R$. That bialgebroid structures on $A$ over $R$ are
in one-to-one correspondence with the strict monoidal forgetful functors $U_A$
was first pointed out by Schauenburg \cite{Schauenburg: Bial}, see also
\cite{Sz}. Of course, the notion of bimonad is much more general than 
bialgebroids. But there is a simple criterion for a bimonad $\mathtt{T}$ on
$_R\M_R$ to be a bialgebroid: The underlying functor $T$ must have a right
adjoint (see Theorem \ref{thm: bgd bimonads}).

Tannaka duality for bialgebroids has recently been proven by Ph\`ung H\^o H\'ai
\cite{Phung Ho Hai} following the tradition of Saavedra \cite{Saavedra},
Deligne \cite{Deligne} and generalizing the results of Ulbrich \cite{Ulbrich},
Schauenburg \cite{Schauenburg: Thesis} for Hopf algebras and of Hayashi
\cite{Hayashi: Tannaka} for weak Hopf algebras. For more about this theory we
refer to \cite{Pareigis: NCAG} and \cite{Joyal-Street: Tannaka} and the
references therein. The approach of the present paper does not fit into this
series mathematically but perhaps `physically'. The categories we think about
are module categories and therefore are not small. This forbids to compute the
(quantum) group(oid) object as a coend. Also, we do not use any finiteness
condition on the images of the functors. Instead we assume that the functors
have left adjoints, and at the end also right ones. The question of when the
bialgebroid we construct has an antipode, in either Lu's \cite{Lu} or
Schauenburg's \cite{Schauenburg: Duals} sense, is not addressed in this paper.

We use the following terminology. A \textit{lax monoidal functor} is a triple
$\bra F,\tau,\iota\ket$ where $F$ is a functor between monoidal categories
from $\bra\C,\oS,\iS\ket$ to $\bra \M,\oT,\iT\ket$,
$\tau_{c,d}\colon Fc\oT Fd\to F(c\oS d)$ is a natural transformation and
$\iota:\iT\to F\iS$ is an arrow such that the usual hexagonal and the two
square diagrams commute \cite{MacLane}. A \textit{lax comonoidal functor} $\bra
F,\sigma,\upsilon\ket$ from $\C$ to $\M$ is a lax monoidal functor from $\C\op$
to $\M\op$ (no change in the monoidal product). That is to say,
$\sigma_{c,d}\colon F(c\oS d)\to Fc\oT Fd$ and $\upsilon\colon F\iS\to\iT$.
The commutative diagrams they satisfy are obtained by reversing all the arrows
is the usual hexagon and square diagrams of a lax monoidal functor. A
\textit{monoidal functor} is a lax monoidal functor in which the $\tau_{c,d}$'s
and the $\iota$ are isomorphisms. \textit{Comonoidal functors} are defined
analogously. Of course, monoidal and comonoidal functors are essentially the
same: $\bra F,\tau,\iota\ket$ is monoidal iff $\bra
F,\tau^{-1},\iota^{-1}\ket$ is comonoidal. 
In later sections where there will be a shortage of the Greek alphabet we
shall use the notation $\bra F,F_2,F_0\ket$ for a lax monoidal and $\bra
F,F^2,F^0\ket$ for a lax comonoidal functor. 
There is also a dual analogue of monoidal natural transformations.
A monoidal natural transformation between lax comonoidal functors will be
called \textit{comonoidal}. So the monoidal categories, the lax comonoidal
functors and the comonoidal natural transformations form a 2-category
$\ComonCat$ just like the lax monoidal functors are the 1-cells of the
2-category $\MonCat$.

A bimonad on the monoidal category $\bra\M,\o,i\ket$ is a quintuple
$\bra T,\gamma,\pi,\mu,\eta\ket$ where $\bra
T,\mu\colon T^2\to T,\eta\colon \M\to T\ket$ is a monad on $\M$,
$\bra T,\gamma\colon T\o\to\o(T\x T),\pi\colon Ti\to i\ket$ is a
comonoidal functor and $\mu$ and $\eta$ satisfy compatibility conditions with
$\gamma$ and $\pi$. The 10 commutative diagrams these natural
transformations should satisfy are the simplest possible ones in spite of the
fact that they are equivalent, at least if $T$ has a right adjoint, to the
somewhat unpleasent bialgebroid axioms that either mix algebra and coalgebra
structures in a painful way \cite{Sz} or use a non-monoidal product in the
definition of coalgebras \cite{Takeuchi}. 

At last but not least the bimonad description offers natural ways to define
the category of bialgebroids which, even for the special case of weak
bialgebras and weak Hopf algebras \cite{BSz,BNSz,Nill}, has not been
investigated in detail yet. Unfortunately, it offers two natural ways. At
first, since bimonads are monads, one can take Street's definition
\cite{Street: Mnd} of the 2-category $\Mnd(\ComonCat)$ of monads in
$\ComonCat$. In this framework all functors are lax comonoidal so the monad
morphisms $\bra G,\varphi\ket$ involve lax comonoidal functors $G$.
This approach was chosen by McCrudden in the preprint \cite{McCr} many
results of which overlap ours. But the present paper uses another
method to define the category of bimonads and therefore of bialgebroids. We
insist on having lax monoidal functors in monad morphisms since it is more
motivated by previous experience with quantum groupoids. So our 2-category
$\Bmd$ of bimonads has objects the bimonads but has lax monoidal functors in
the definition of monad morphisms. This choice forces the $\varphi$ in the
monad morphism to be \textit{ambimonoidal} natural transformation (see
Definition \ref{def: ambi}). With this tentative expression we refer to the
unique way of compatibility with the (co)monoidal structures which, however,
lies between the usual monoidality and comonoidality of natural
transformations.

Having a 2-category of bialgebroids one can imagine the 2-category of Hopf
algebroids as a subcategory. The equivalence classes of objects in 
this 2-category may turn out to be the appropriate objects
which characterize a class of monoidal categories uniquely, similarly to the
Doplicher-Roberts characterization \cite{DR} of certain symmetric monoidal
C$^*$-categories as representation categories of uniquely determined compact
groups.

\section{Monoidal adjunctions}

\subsection{Lax monoidal functors with left adjoint}

The whole content of this paper rests on the following result
which can be found in various forms in the categorical literature,
in most general form perhaps in \cite{Kelly}. Nevertheless we provide an
explicit proof in order for the paper to be selfcontained even for
non-specialists.

\begin{thm} \label{thm: F-sigma-pi}
Let $\bra U,\tau,\iota\ket$ be a lax monoidal functor from $\bra
\C,\oS,\iS\ket$ to $\bra\M,\oT,\iT\ket$. If $U$ has a left adjoint
$F\colon\M\to\C$ with counit $\eps\colon FU\to\C$ and unit $\eta\colon \M\to
UF$ then the natural transformation
\beq
\sigma_{x,y}:=\eps_{Fx\oS Fy}\circ F\tau_{Fx,Fy}\circ F(\eta_x\oT \eta_y)
\eeq
and the arrow 
\beq
\upsilon:=\eps_\iS\circ F\iota
\eeq
give rise to a lax comonoidal functor $\bra F,\sigma,\upsilon\ket$.
\end{thm}
\begin{proof}
The expression for $\sigma$ contains three products: the horizontal and
vertical composition of natural transformations and the Cartesian product. So
the statement requires 3-categorical computations, in the monoidal 2-category
$\Cat$. Denoting the vertical composition by $\bullet$, the horizontal by
juxtaposition and the Cartesian product by $\x$ and using the, perhaps strange,
precedence of horizontal composition coming first and vertical coming last, we
can write
\beq
\sigma=\eps\oS(F\x F)\bullet F\tau(F\x F)\bullet F\oT(\eta\x\eta)\ .
\eeq
Denoting by $\a$, $\l$, $\r$ the associativity and unit coherence isomorphisms
in any one of the monoidal categories, the proof of the hexagon
\beq
\a(F\x F\x F)\bullet\oS(F\x \sigma)\bullet\sigma(\M\x\oT)
=\oS(\sigma\x F)\bullet\sigma(\oT\x\M)\bullet F\a
\eeq
goes as follows.
\begin{align}
\begin{split}
& \oS(\sigma\x F)\bullet\sigma(\oT\x\M)=\\
=&\oS(\eps\oS(F\x F)\x F)\bullet\oS(F\tau(F\x F)\x
F)\bullet\oS(F\oT(\eta\x\eta)\x F)\\
&\bullet\eps\oS(F\x F)(\oT\x\M)\bullet F\tau(F\x F)(\oT\x\M)\bullet
F\oT(\eta\x\eta)(\oT\x\M)\\
=&\eps\oS(\oS\x\C)(F\x F\x F)\bullet FU\oS(\eps\oS(F\x F)\x F)\bullet
FU\oS(F\tau(F\x F)\x F)\\
&\bullet FU\oS(F\oT(\eta\x\eta)\x F)\bullet F\tau(F\oT\x F)\bullet
F\oT(\eta\oT\x\eta)\\
=&\eps\oS(\oS\x\C)(F\x F\x F)\bullet FU\oS(\eps\oS(F\x F)\x F)\bullet
FU\oS(F\tau(F\x F)\x F)\\
&\bullet F\tau(F\oT(UF\x UF)\x F)\bullet F\oT(UF\oT(\eta\x\eta)\x UF)\bullet
F\oT(\eta\oT\x\eta)\\
=&\eps\oS(\oS\x\C)(F\x F\x F)\bullet FU\oS(\eps\oS(F\x F)\x F)\bullet
FU\oS(F\tau(F\x F)\x F)\\
&\bullet F\tau(F\oT(UF\x UF)\x F)\bullet F\oT(\eta\oT(UF\x UF)\x UF) \bullet
F\oT(\oT(\eta\x\eta)\x\eta)\\
=&\eps\oS(\oS\x\C)(F\x F\x F)\bullet FU\oS(\eps\oS(F\x F)\x F)\bullet
F\tau(FU\oS(F\x F)\x F)\\
&\bullet F\oT(UF\tau(F\x F)\x UF)\bullet F\oT(\eta\oT(UF\x FF)\x UF) \bullet
F\oT(\oT(\eta\x\eta)\x\eta)\\
\end{split}
\end{align}
\begin{align*}
\begin{split}
=&\eps\oS(\oS\x\C)(F\x F\x F)\bullet F\tau(\oS(F\x F)\x F)\bullet
F\oT(U\eps\oS(F\x F)\x UF)\\
&\bullet F\oT(UF\tau(F\x F)\x UF)\bullet F\oT(\eta\oT(UF\x UF)\x UF) \bullet
F\oT(\oT(\eta\x\eta)\x\eta)\\
=&\eps\oS(\oS\x\C)(F\x F\x F)\bullet F\tau(\oS(F\x F)\x F)\bullet
F\oT(U\eps\oS(F\x F)\x UF)\\
&\bullet F\oT(\eta U\oS(F\x F)\x UF)\bullet F\oT(\tau(F\x F)\x UF) \bullet
F\oT(\oT(\eta\x\eta)\x\eta)\\
=&\eps\oS(\oS\x\C)(F\x F\x F)\bullet F\tau(\oS\x\C)(F\x F\x F)\\
&\bullet F\oT(\tau\x U)(F\x F\x F) \bullet F\oT(\oT(\eta\x\eta)\x\eta)
\end{split}
\end{align*}
where in the subsequent equations we used the definition of $\sigma$, 
naturality of the $\eps$ of the 4th term, naturality of the $\tau$ of the 5th
term,  naturality of the first $\eta$ of the 6th term, naturality of the
$\tau$ of the 4th term, naturality of the $\eta$ of the 5th term, and at last
the adjointness $F\dashv U$ and some cosmetics.
A similarly long calculation, or arguing with $\x^{op}$, yields the formula
\begin{align}
\begin{split}
&\oS(F\x \sigma)\bullet\sigma(\M\x\oT)=\\ 
=&\eps\oS(\C\x\oS)(F\x F\x F)\bullet F\tau(\C\x\oS)(F\x F\x F)\\
&\bullet F\oT(U\x \tau)(F\x F\x F)\bullet F\oT(\eta\x\oT(\eta\x\eta))
\end{split}
\end{align}
Composing the first with $F\a$ and the second with $\a(F\x F\x F)$ and using
the hexagon for $\tau$
\beq
\tau(\oS\x\C)\bullet\oT(\tau\x U)\bullet\a(U\x U\x U)
=U\a\bullet\tau(\C\x\oS)\bullet\oT(U\x\tau)
\eeq
one immediately obtains the hexagon for $\sigma$.

It remained to show the squares of $\sigma$ and $\upsilon$
\bea
\l_{Fx}\c (\upsilon\oS Fx)\c \sigma_{\iT,x}&=&F\l_x\\
\r_{Fx}\c (Fx\oS \upsilon)\c \sigma_{x,\iT}&=&F\r_x
\eea
 but we will suffice with proving the first. 
\begin{align}
\begin{split}
&\l_{Fx}\c (\upsilon\oS Fx)\c \sigma_{\iT,x}=\\
&=\l_{Fx}\c(\eps_\iS\oS Fx)\c(F\iota\oS Fx)\c\eps_{F\iT\oS Fx}\c
F\tau_{F\iT,Fx} \c F(\eta_\iT\oT\eta_x)\\
&=\l_{Fx}\c(\eps_\iS\oS Fx)\c\eps_{FU\iS\oS Fx}\c FU(F\iota\oS Fx)\c
F\tau_{F\iT,Fx} \c F(\eta_\iT\oT\eta_x)\\
&=\l_{Fx}\c\eps_{\iS\oS Fx}\c FU(\eps_\iS\oS Fx)\c F\tau_{FU\iS,Fx}\c
F(UF\iota\oT UFx)\c F(\eta_\iT\oT\eta_x)\\
&=\eps_{Fx}\c FU\l_{Fx}\c F\tau_{\iS,Fx}\c F(U\eps_\iS\oT UFx)\c
F(\eta_{U\iS}\oT UFx)\c F(\iT\oT\eta_x)\\
&=\eps_{Fx}\c F\left(U\l_{Fx}\c\tau_{\iS,Fx}\c (\iota\oT UFx)\right)\c
F(\iT\oT\eta_x)\\
&=\eps_{Fx}\c F\eta_x\c F\l_{Fx}\ =\ F\l_x
\end{split}
\end{align}
where in the last but one equality the square of $\tau$ and $\iota$ was
used.
\end{proof}

We remark that the converse of the above Theorem holds, too. If $\bra
F,\sigma,\upsilon\ket$ is lax comonoidal then 
\begin{align}
\tau_{a,b}&:=U(\eps_a\oS\eps_b)\c U\sigma_{Ua,Ub}\c\eta_{Ua\oT Ub}\\
\iota&:=U\upsilon\c\eta_\iV
\end{align}
defines a lax monoidal structure for $U$. However, monoidality of $U$, i.e.,
invertiblity of $\tau$ and $\iota$ does not imply comonoidality of $F$, i.e.,
invertibility of $\sigma$ and $\upsilon$.

\subsection{Comonoidality of the unit and counit}

We specialize the above Theorem to monoidal $U$. Then we get the following.

\begin{pro} \label{pro: comon adj}
If $\bra U,\tau,\iota\ket$ is monoidal and $F$ is a left adjoint of $U$ then
the $\bra F,\sigma,\upsilon\ket$ given in Theorem \ref{thm: F-sigma-pi} is
the unique lax comonoidal structure on $F$ for which the given adjunction data
$\eps$ and $\eta$ are comonoidal natural transformations, i.e., for which
\begin{align} \label{dia: comon eps}
\begin{CD}
FU(a\oS b)@>\sigma_{Ua,Ub}\c F\tau^{-1}_{a,b}>>FUa\oS FUb\\
@V{\eps_{a\oS b}}VV  @VV{\eps_a\oS \eps_b}V\\
a\oS b @=a\oS b
\end{CD}
&\qquad
\begin{CD}
FU\iS@>F\iota^{-1}>>F\iT\\
@V{\eps_\iS}VV @VV{\upsilon}V\\
\iS@= \iS
\end{CD}
\\
\intertext{and} \label{dia: comon eta}
\begin{CD}
x\oT y @= x\oT y\\
@V{\eta_{x\oT y}}VV  @VV{\eta_x\oT\eta_y}V\\
UF(x\oT y)@>\tau^{-1}_{Fx,Fy}\c U\sigma_{x,y}>>UFx\oT UFy
\end{CD}
&\qquad
\begin{CD}
\iT@= \iT\\
@V{\eta_\iT}VV @AA{\iota^{-1}}A\\
UF\iT@>U\upsilon>>U\iS
\end{CD}
\end{align}
are commutative.
\end{pro}
\begin{proof}
Using invertibility of $\tau_{a,b}$ and $\iota$ the above diagrams can be read
as the equations
\begin{align} \label{eq: monadj1}
(\eps_a\oS\eps_b)\circ\sigma_{Ua,Ub}&=\eps_{a\oS b}\circ F\tau_{a,b}\\
\upsilon&=\eps_\iS\circ F\iota \label{eq: monadj2}\\
\tau_{Fx,Fy}\circ(\eta_x\oT\eta_y)&=U\sigma_{x,y}\circ \eta_{x\oT y}
\label{eq: monadj3}\\
\iota&=U\upsilon\circ\eta_\iT\label{eq: monadj4}
\end{align}
where $a,b$ run over the objects of $\C$ and $x,y$ over those of $\M$.
Now we are left to show that these equations have a unique solution for
$\sigma$ and $\upsilon$. For $\upsilon$ this is obvious from the second
equation. In order to obtain $\sigma$ apply $F$ to the third equation and
multiply the result with $\eps_{Fx,Fy}$.
\begin{align*}
\eps_{Fx,Fy}\circ F\tau_{Fx,Fy}\circ F(\eta_x\oT\eta_y)&=\eps_{Fx,Fy}\circ
FU\sigma_{x,y}\circ F\eta_{x\oT y}\\
&=\sigma_{x,y}\circ\eps_{F(x\oT y)}\circ F\eta_{x\oT y}=\sigma_{x,y}
\end{align*}
which is indeed the comonoidal structure of Theorem \ref{thm: F-sigma-pi}.
\end{proof}

In the course of the above proof we have seen that equations (\ref{eq:
monadj1}-\ref{eq: monadj4}) have unique solutions for $\sigma$ and $\upsilon$ if
the $\eps$, $\eta$, $\tau$, $\iota$ are given. This is true even if $U$ is lax
monoidal although the comonoidality diagrams (\ref{dia: comon eps}) and
(\ref{dia: comon eta}) loose their meaning. (The $\eps$ and $\eta$ are
ambimonoidal, however, in the sense of Definition \ref{def: ambi}.) Similarly,
equations (\ref{eq: monadj1}-\ref{eq: monadj4}) can be solved for $\tau$ and
$\iota$ if the others are given. This motivates the

\begin{defi}
A lax monoidal functor $U$ and a lax comonoidal functor $F$ is
called a monoidal adjoint pair, or simply a monoidal adjunction $F\dashv U$, if
$F$ is the left adjoint of $U$ as ordinary functors and the counit $\eps\colon
FU\to \C$ and the unit $\eta\colon\M\to UF$ can be chosen to satisfy equations
(\ref{eq: monadj1}), (\ref{eq: monadj2}), (\ref{eq: monadj3}) and (\ref{eq:
monadj4}). 
\end{defi}

\begin{lem} \label{lem: uniq mon adj}
If $U$ is monoidal and $F\dashv U$, $F'\dashv U$ are two monoidal adjunctions
then there exists a comonoidal natural isomorphism $F\iso F'$.
\end{lem}
\begin{proof}
As in the case of the proof of uniqueness of left adjoints of functors up to
natural isomorphisms one takes the natural isomorphism $\eps F'\b F\eta'\colon
F\to F'$ which is made of horizontal and vertical composites of lax comonoidal
functors and of comonoidal natural transformations, by Proposition \ref{pro:
comon adj}, so it is itself comonoidal.
\end{proof}

\section{The monoidal Eilenberg-Moore construction}

\subsection{Monoidal functors with left adjoints and bimonads}

The continuing assumption is that $\bra U, \tau,\iota\ket$ is a monoidal
functor from $\bra \C,\oS,\iS\ket$ to $\bra\M,\oT,\iT\ket$ and that the
functor $U$ has a left adjoint $F\colon\M\to\C$ with counit $\eps\colon FU\to
\C$ and unit $\eta\colon \M\to UF$. In this situation Proposition \ref{pro:
comon adj} has the following 

\begin{cor} \label{cor: mon adj}
If $\bra U,\tau,\iota\ket$ is a monoidal functor with the underlying functor
$U$ having a left adjoint then the monad $\bra T,\mu,\eta\ket$ associated to
the adjunction data $\eps\colon FU\to\C$, $\eta\colon\M\to UF$ is such that
$T=UF$ is a lax comonoidal functor and $\mu\colon T^2\to T$ and $\eta\colon\M
\to T$ are comonoidal natural transformations.
\end{cor}
\begin{proof}
By Theorem \ref{thm: F-sigma-pi} the left adjoint $F\colon\M\to\C$ has a
comonoidal structure. $U$ being monoidal it is also comonoidal and the
composition $T=UF$ of lax comonoidal functors is again lax comonoidal. In
Proposition \ref{pro: comon adj} we have seen that $\eps$ and $\eta$ are
comonoidal natural transformations w.r.t. this comonoidal structure on $U$ and
$F$. Therefore $\mu=U\eps F$ is also comonoidal since the monoidal categories,
the lax comonoidal functors and the comonoidal natural transformations form a
2-category.
\end{proof}

The monad we have obtained in the above Corollary suggests the following
\begin{defi} \label{def: bmd}
Let $\M$ be a monoidal category. Then a monoid $\bra T,\mu,\eta\ket$ in the
monoidal category of lax comonoidal endofunctors $\M\to\M$ is called a
bimonad in $\M$. 
\end{defi}
 
That is to say, a bimonad consists of 6 items,
\begin{enumerate}
\item a monoidal category $\bra \M,\oT,\iT\ket$
\item a functor $T\colon\M\to\M$
\item a natural transformation $\gamma_{x,y}\colon T(x\oT y)\to Tx\oT Ty$
\item an arrow $\pi\colon T\iT\to\iT$
\item a natural transformation $\mu_x\colon T^2x\to Tx$
\item and a natural transformation $\eta_x\colon x\to Tx$
\end{enumerate}
subjected to satisfy 6 axioms in the form of 10 commutative diagrams:
\begin{description}
\item[BMD 1] $\gamma$ is coassociative,
\beq \label{dia: BMD 1}
\begin{CD}
T(x\oT (y\oT z))@>\gamma_{x,y\oT z}>>Tx\oT T(y\oT z)@>Tx\oT\gamma_{y,z}>>
Tx\oT(Ty\oT Tz)\\
@V{T\a_{x,y,z}}VV @. @VV{\a_{Tx,Ty,Tz}}V\\
T((x\oT y)\oT z)@>\gamma_{x\oT y,z}>>T(x\oT y)\oT Tz@>\gamma_{x,y}\oT Tz>>
(Tx\oT Ty)\oT Tz
\end{CD}
\eeq
\item[BMD 2] $\pi$ is a counit for $\gamma$,
\beq \label{dia: BMD 2}
\begin{CD}
T(\iT\oT x)@>\gamma_{\iT,x}>>T\iT\oT Tx\\
@V{\wr}VV @VV{\pi\oT Tx}V\\
Tx@<\sim<<\iT\oT Tx
\end{CD}
\qquad\quad\text{and}\qquad
\begin{CD}
T(x\oT\iT)@>\gamma_{x,\iT}>>Tx\oT T\iT\\
@V{\wr}VV @VV{Tx\oT\pi}V\\
Tx@<\sim<<Tx\oT\iT
\end{CD}
\eeq
\item[BMD 3] $\mu$ is comonoidal,
\beq \label{dia: BMD 3}
\begin{CD}
T^2(x\oT y)@>\gamma_{Tx,Ty}\c T\gamma_{x,y}>>T^2x\oT T^2y\\
@V{\mu_{x\oT y}}VV  @VV{\mu_x\oT\mu_y}V\\
T(x\oT y)@>\gamma_{x,y}>>Tx\oT Ty
\end{CD}
\qquad\quad\text{and}\qquad
\begin{CD}
T^2\iT@>\pi\c T\pi>>\iT\\
@V{\mu_\iT}VV  @|\\
T\iT@>\pi>>\iT
\end{CD}
\eeq
\item[BMD 4] $\eta$ is comonoidal,
\beq \label{dia: BMD 4}
\begin{CD}
x\oT y@= x\oT y\\
@V{\eta_{x\oT y}}VV  @VV{\eta_x\oT\eta_y}V\\
T(x\oT y)@>\gamma_{x,y}>> Tx\oT Ty
\end{CD}
\qquad\quad\text{and}\qquad
\begin{CD}
\iT@=\iT\\
@V{\eta_\iT}VV @|\\
T\iT@>\pi>>\iT
\end{CD}
\eeq
\item[BMD 5] $\mu$ is associative,
\beq \label{dia: BMD 5}
\begin{CD}
T^3x@>T\mu_x>>T^2x\\
@V{\mu_{Tx}}VV @VV{\mu_x}V\\
T^2x@>\mu_x>>Tx
\end{CD}
\eeq
\item[BMD 6] $\eta$ is a unit for $\mu$,
\beq \label{dia: BMD 6}
\begin{CD}
Tx@>\eta_{Tx}>>T^2x\\
@| @V{\mu_x}VV\\
Tx@=Tx
\end{CD}
\qquad\quad\text{and}\qquad
\begin{CD}
Tx@>T\eta_x>>T^2x\\
@| @V{\mu_x}VV\\
Tx@= Tx
\end{CD}
\eeq
\end{description}
 
Thus Corollary tells us that every monoidal adjunction $F\dashv U$,
with $U$ monoidal, determines a bimonad with underlying monad the classical
construction $T=\bra UF,U\eps F,\eta\ket$. Explicitely, if $\bra
U,\tau,\iota\ket$ is a monoidal functor and $F\dashv U$ is an ordinary
adjunction with unit $\eta$ and counit $\eps$ then the associated bimonad is
this.  
\begin{align}
T&:=UF\\ 
\gamma_{x,y}&:=\tau^{-1}_{Fx,Fy}\c U\eps_{Fx\oS Fy}\c
UF\tau_{Fx,Fy}\c UF(\eta_x\oT \eta_y) \label{eq: gamma}\\  
\pi&:=\iota^{-1}\c U\eps_\iS\c UF\iota \label{eq: pi}\\
\mu_x&:=U\eps_{Fx}\\
\eta_x&:=\eta_x\,.
\end{align}

\subsection{The monoidal Eilenberg-Moore category}

In this subsection $\M$ is a monoidal category and $\bra T,\mu,\eta\ket$ is a
bimonad on $\M$.

The Eilenberg-Moore category $\M^T$ has as objects the $T$-algebras, i.e., pairs
$\bra x,\alpha\ket$ where $x$ is an object in $\M$ and $\alpha\colon Tx\to x$
satisfies
\beq \label{eq: T-alg}
\begin{CD}
T^2x@>T\alpha>> Tx\\
@V{\mu_x}VV @VV{\alpha}V\\
Tx@>\alpha>> x
\end{CD}
\qquad\qquad
\begin{CD}
x@>\eta_x>> Tx\\
@| @VV{\alpha}V\\
x@=x
\end{CD}
\eeq
The arrows from $\bra x,\alpha\ket$ to $\bra y,\beta\ket$ are the arrows
$t\colon x\to y$ in $\M$ such that
\beq
\begin{CD}
Tx@>Tt>>Ty\\
@V{\alpha}VV @VV{\beta}V\\
x@>t>>y
\end{CD}
\eeq
is commutative. The functor 
\beq
U^T\colon\M^T\to \M\,,\qquad \bra x,\alpha\ket\mapsto x
\eeq
is called the Eilenberg-Moore forgetful functor.

\begin{pro} \label{pro: mon EM}
Let $\bra T,\gamma,\pi,\mu,\eta\ket$ be a bimonad. Then its
Eilenberg-Moore category $\M^T$ has the following monoidal structure. For
$T$-algebras $\bra x,\alpha\ket$ and $\bra y,\beta\ket$ let their tensor
product be
\beq \label{eq: oEM}
\bra x,\alpha\ket\oEM\bra y,\beta\ket:=\bra x\oT
y,(\alpha\oT\beta)\circ\gamma_{x,y}\ket
\eeq
The tensor product of $T$-algebra arrows coincides with their tensor product
as arrows in $\M$. Then $\oEM$ gives rise to a monoidal structure on $\M^T$
such that the forgetful functor $U^T\colon\M^T\to\M$ becomes strictly monoidal.
\end{pro}
\begin{proof}
In order to show that (\ref{eq: oEM}) is really a $T$-algebra we need to verify
the two defining diagrams of (\ref{eq: T-alg}). The first of these follow from
\begin{align*}
&(\alpha\oT\beta)\circ\gamma_{x,y}\circ T(\alpha\oT\beta)\circ T\gamma_{x,y}=\\
=&(\alpha\oT\beta)\circ(T\alpha\oT T\beta)\circ\gamma_{Tx,Ty}\circ T\gamma_{x,y}=\\
=&(\alpha\oT\beta)\circ(\mu_x\oT\mu_y)\circ\gamma_{Tx,Ty}\circ T\gamma_{x,y}=\\
=&(\alpha\oT\beta)\circ\gamma_{x,y}\circ\mu_{x\oT y}
\end{align*}
where in the last equation we used comonoidality of $\mu$. The second diagram
follows from
\begin{align*}
&(\alpha\oT\beta)\circ\gamma_{x,y}\circ\eta_{x\oT y}=\\
=&(\alpha\oT\beta)\circ(\eta_x\oT\eta_y)=x\oT y
\end{align*}
where comonoidality of $\eta$ had to be used.

For $T$-algebra arrows $t\colon\bra x,\alpha\ket\to\bra x',\alpha'\ket$ and
$s\colon\bra y,\beta\ket\to\bra y',\beta'\ket$ we defined
\beq
t\oEM s:=t\oT s
\eeq
which is indeed a $T$-algebra arrow because
\begin{align*}
&(\alpha'\oT\beta')\circ\gamma_{x',y'}\circ T(t\oT s)=\\
=&(\alpha'\oT\beta')\circ(Tt\oT Ts)\circ \gamma_{x,y}=\\
=&(t\oT s)\circ(\alpha\oT\beta)\circ \gamma_{x,y}
\end{align*}
This finishes the definition of the functor $\oEM$. As for the monoidal unit
we set 
\beq
\iT^T:=\bra \iT,\pi\ket
\eeq
Now we are going to show that the coherence isomorphism $\a$, $\l$, and $\r$
of $\bra\M,\oT,\iT\ket$, when considered as arrows in $\M^T$, serve as 
coherence isomorphisms of $\bra\M^T,\oEM,\iT^T\ket$. For three
$T$-algebras $x^T=\bra x,\alpha_x\ket$, $y^T=\bra y,\alpha_y\ket$ and
$z^T=\bra z,\alpha_z\ket$ we have
\begin{align*}
x^T\oTT (y^T\oTT z^T)&=\bra x\oT (y\oT z),
(\alpha_x\oT(\alpha_y\oT\alpha_z))\circ (Tx\oT\gamma_{y,z})\circ
\gamma_{x,y\oT z}\ket\\
(x^T\oTT y^T)\oTT z^T&=\bra (x\oT y)\oT z,
((\alpha_x\oT\alpha_y)\oT\alpha_z)\circ (\gamma_{x,y}\oT Tz)\circ
\gamma_{x\oT y,z}\ket
\end{align*}
and the calculation
\begin{align*}
&\a_{x,y,z}\circ(\alpha_x\oT(\alpha_y\oT\alpha_z))\circ
(Tx\oT\gamma_{y,z})\circ \gamma_{x,y\oT z}=\\
=&((\alpha_x\oT\alpha_y)\oT\alpha_z)\circ\a_{Tx,Ty,Tz}\circ
(Tx\oT\gamma_{y,z})\circ \gamma_{x,y\oT z}=\\
=&((\alpha_x\oT\alpha_y)\oT\alpha_z)\circ(\gamma_{x,y}\oT Tz)\circ
\gamma_{x\oT y,z}\circ T\a_{x,y,z}
\end{align*}
proves that $\a_{x,y,z}$ is an isomorphism
\beq
x^T\oTT (y^T\oTT z^T)\rarr{\a_{x,y,z}}(x^T\oTT y^T)\oTT z^T
\eeq
of $\M^T$, indeed. In order to show that $\iT^T$ is a left unit notice that
\beq
\iT^T\oTT x^T=\bra \iT\oT x,(\pi\oT\alpha_x)\circ\gamma_{\iT,x}\ket
\eeq
and therefore
\begin{align*}
&\l_x\circ(\pi\oT\alpha_x)\circ\gamma_{\iT,x}=\\
=&\l_x\circ(\iT\oT\alpha_x)\circ(\pi\oT Tx)\circ\gamma_{\iT,x}=\\
=&\l_{Tx}\circ(\pi\oT Tx)\circ\gamma_{\iT,x}=\\
=&T\l_x\circ\alpha_x
\end{align*}
proves that $\l_x$ is an isomorphism
\beq
\iT^T\oTT x^T\rarr{\l_x}x^T
\eeq
of $\M^T$. Similarly, $\r_x\colon x^T\oTT\iT^T\iso x^T$ for all objects $x^T$
of $\M^T$. This finishes the construction of a monoidal structure
$\bra\M^T,\oTT,\iT^T\ket$ on the Eilenberg-Moore category.
It is clear from the construction that the forgetful functor $U^T$ is
strictly monoidal.
\end{proof}

\begin{pro} \label{pro: F^T}
Let $T$ be a bimonad. Then the strict monoidal $U^\ttT$ has a lax comonoidal
left adjoint, the free $T$-algebra functor
\beq
F^\ttT\colon\M\to\M^\ttT\,,\quad x\mapsto\bra Tx,T^2x\rarr{\mu_x}Tx\ket\,.
\eeq
such that $U^\ttT F^\ttT=T$, as lax comonoidal functors.
\end{pro}
\begin{proof}
Left adjointness is proven as in the textbooks. As for the lax comonoidal
structure notice that $\gamma_{x,y}$ provides a $T$-algebra arrow from
$F^\ttT(x\o y)$ to
\[
F^\ttT x\oalg F^\ttT y=\bra Tx\o Ty,(\mu_x\o\mu_y)\c\gamma_{Tx,Ty}\ket
\]
because $\gamma_{x,y}\c\mu_{x\o y}=(\mu_x\o\mu_y)\c\gamma_{Tx,Ty}\c
T\gamma_{x,y}$ which is precisely the first diagram in (\ref{dia: BMD 3}).
Similarly, the second diagram of (\ref{dia: BMD 3}) is the condition
for  
\[
\hat\pi\colon\ F^\ttT=\bra Ti,\mu_i\ket\rarr{\pi}\bra i,\pi\ket=i^\ttT
\]
to be a $T$-algebra arrow. Now (\ref{dia: BMD 1}-\ref{dia: BMD 2}) imply that
the triple $\bra F^\ttT,\hat\gamma,\hat\pi\ket$ is a lax comonoidal functor.
The relations $U^\ttT F^\ttT=T$, $U^\ttT\hat\gamma=\gamma$ and
$U^\ttT\hat\pi=\pi$ are obvious. 
\end{proof}

\subsection{The Tannakian reconstruction}

Recall that the Eilenberg-Moore comparison functor $K\colon \C\to\M^T$ maps
the objects $c$ of $\C$ into the $T$-algebras
\beq
Kc\ :=\ \bra Uc, U\eps_c\colon TUc\to Uc\ket
\eeq
and the arrows $\psi\colon c\to d$ to the $T$-algebra morphisms 
\beq
K\psi\ :=\ \bra Uc,U\eps_c\ket\rarr{U\psi}\bra Ud,U\eps_d\ket\,.
\eeq
This functor allows to factorize the given $U$ through the category of
$T$-algebras as $U^TK=U$.

\begin{pro} \label{pro: K is mon}
Let $U\colon\C\to\M$ be a monoidal functor with left adjoint and $T$ be the
associated bimonad. Then the Eilenberg-Moore comparison functor
$K\colon\C\to\M^T$ has a unique (lax) monoidal structure such that the
factorization $U=U^TK$ is a factorization of monoidal functors.
\end{pro}
\begin{proof}
Since $U^T$ is strict monoidal, if $\bra K,\hat\tau,\hat\iota\ket$ is a lax
monoidal functor such that $U=U^TK$ then $U\hat\tau=\tau$ and
$U\hat\iota=\iota$. That is to say, the unique lax monoidal structure, if
exists, is monoidal and it is obtained by lifting the arrows $\tau_{c,d}$ and
$\iota$ to $T$-algebra arrows. Taking into account formula (\ref{eq: gamma})
the action in the tensor product
\[
Kc\oalg Kd\ =\ \bra Uc\oT Ud,(U\eps_c\oT U\eps_d)\c\gamma_{Uc,Ud}\ket
\]
can be writen as
\begin{align*}
(U\eps_c\oT U\eps_d)\c\gamma_{Uc,Ud}&=\tau_{c,d}^{-1}\c U(\eps_c\oS\eps_d)\c
U\sigma_{Uc,Ud}\\
&=\tau_{c,d}^{-1}\c U\eps_{c\oS d}\c T\tau_{c,d}
\end{align*}
where in the last equation the monoidal adjunction (\ref{eq: monadj1}) has
been used. This result, up to multiplying with $\tau_{c,d}$, is precisely the
lifting condition for $\tau_{c,d}$ to be $T$-algebra morphism $Kc\oalg Kd\to
K(c\oS d)$. As for the unit map $\iota\colon \iT\to U\iS$, it has a lift to a
$T$-algebra morphism
\[
\bra\iT,\pi\ket\to\bra U\iS,U\eps_\iS\ket\ =\ K\iS
\]
if and only if $\iota\c\pi=U\eps_\iS\c T\iota$. The right hand side is equal
to $U\upsilon$ by (\ref{eq: monadj2}) and the left hand side is equal to
$U\upsilon$ by (\ref{eq: pi}).
\end{proof}

Our main theorem of the Tannakian type relates adjointable monoidal functors
to bimonads in the following way.
\begin{thm} \label{thm: 1-univ}
Let $U\colon \C\to\M$ be a monoidal functor possessing a left adjoint.
Then there exists a bimonad $\mathbf{T}$ on $\M$
and a monoidal functor $K\colon\C\to\M^\ttT$ such that 
\begin{enumerate}
\item $U$ has the monoidal factorization $U=U^\ttT K$ where
$U^\ttT\colon\M^\ttT\to\M$ is the monoidal Eilenberg-Moore forgetful functor,
\item the pair $\bra T,K\ket$ is universal with respect to property (1).
That is to say, if $\mathbf{S}$ is a bimonad on
$\M$ and $L\colon \A\to\M^\ttS$ is a lax monoidal functor such that $U=U^\ttS
L$, as lax monoidal functors, then there exists a unique natural
transformation $\varphi\colon S\to T$ such that 
\begin{enumerate}
\item $U^\varphi K=L$ where $U^\varphi$ is the functor mapping a
$T$-algebra $\bra x,\alpha\ket$ to the $S$-algebra $\bra x,\alpha\c
\varphi_x\ket$, 
\item $\varphi$ is comonoidal,
\item $\varphi$ is a monad morphism.
\end{enumerate}
\end{enumerate}
\end{thm}
\begin{proof}
As for the existence of $T$ and $K$ with property (1) one takes for $T$
the bimonad associated to $U$ and to one of its left adjoints $F$ by Corollary
\ref{cor: mon adj} and for $K$ the Eilenberg-Moore comparison functor. 
In order to show the universal property (2) we need to prove existence and
uniqueness of $\varphi$. Notice that the functors
$L$ with property $U^\ttS L=U$ can be written as $Lc=\bra Uc,\beta_c\ket$
where functoriality implies that $\beta_c\colon SUc\to Uc$ is natural in
$c\in\C$. Furthermore, lax monoidality of $L$, strict monoidality of $U^\ttS$
and monoidality of the factorization $U=U^\ttS L$ implies that $L$ is monoidal
and $L_2\colon \oalg(L\x L)\iso L\oS$ is the lift of $U_2\colon \oT(U\x U)\iso
U\oS$. Similarly $L_0\colon \bra \iT,S^0\ket\iso L\iS=\bra
U\iS,\beta_\iS\ket$ is the lift of $U_0\colon \iT\iso U\iS$. Since
\begin{align*}
Lc\oalg Ld&=\bra Uc\oT Ud,(\beta_c\oT\beta_d)\c S_{Uc,Ud}\ket\\
L\iS&=\bra U_e,\beta_\iS\ket\,,
\end{align*}
the lifting conditions for $U_{c,d}$ and $U_0$, respectively take the form
\begin{align}
U_{c,d}\c (\beta_c\oT\beta_d)\c S_{Uc,Ud}&=\beta_{c\oS d}\c SU_{c,d}\\
U_0\c S^0&=\beta_\iS\c SU_0
\end{align}
and these are precisely the conditions for $\beta$ to be comonoidal. Let
$\L(U)$ be the category with objects the pairs $\bra S,\beta\ket$ where
$S\colon\M\to\M$ is a lax comonoidal functor $\bra S,S^2,S^0\ket$ and
$\beta\colon SU\to U$ is a comonoidal natural transformation. The arrows from
$\bra R,\alpha\ket$ to $\bra S,\beta\ket$ are the comonoidal natural
transformations $\varphi\colon R\to S$ satisfying $\beta\b \varphi U=\alpha$.
Now it is standard universal algebra to show that if $U$ has a left adjoint
then $\L(U)$ has terminal objects. Moreover, in a terminal object $\bra
T,\omega\ket$ the $T$ is a monad and $\omega$ is an action of $T$. If $\bra
S,\beta\ket$ is an object in which $S$ is a monad and $\beta$ is an action of
$S$ on $U$ then the unique arrow $\varphi\colon \bra S,\beta\ket\to\bra
T,\omega\ket$ is a monad morphism, i.e., satisfies
\begin{align*}
\mu^\ttT\b T\varphi\b\varphi S&=\mu^\ttS\\
\varphi\b\eta^\ttS&=\eta^\ttT\,.
\end{align*}
Now it is easy to see that condition (2) is just the expression of the fact
that the pair $\bra T, \omega\ket$, in which $\omega$ correponds to the
comparison functor $K$, is a terminal object in $\L(U)$. We omit the details
because we shall prove in Theorem \ref{thm: 2-univ} a more general universality
property involving bimonads $S$ on any other monoidal category $\N$. (Cf. also
the proof of Lemma \ref{lem: Q0} or the literature \cite{Pareigis: NCAG}.) 
\end{proof}

The next theorem serves as a characterization of the forgetful functors of
bimonads.

\begin{thm} \label{thm: rep}
Let $\C$ and $\M$ be monoidal categories. For a functor $U\colon\C\to\M$ the
following conditions are equivalent: 
\begin{enumerate}
\item There exists a bimonad $T$ on $\M$ and a monoidal equivalence $K\colon
\C\to\M^\ttT$ such that -- via this equivalence -- $U$ is isomorphic to the
forgetful functor $U^T$.   
\item $U$ is monadic and monoidal. 
\end{enumerate}
\end{thm}
\begin{proof}
$(2)\Rightarrow (1)$ Monadicity of $U$ is by definition the requirement that
$U$ has a left adjoint and the comparison functor $K$ is a category
equivalence. A category equivalence is always part of an adjoint equivalence
\cite[Theorem IV. 4. 1]{MacLane} so there exists a right adjoint of $K$ with
invertible unit and counit. Now $K$ is monoidal, hence comonoidal, therefore
the converse of Theorem \ref{thm: F-sigma-pi} provides a lax monoidal
structure on the right adjoint which, by invertibility of the unit and counit,
is actually monoidal. This proves that $K$ is a monoidal equivalence and the
rest, $U=U^\ttT K$, is obvious.

$(1)\Rightarrow (2)$ The Eilenberg-Moore forgetful functor $U^\ttT$ is always
monadic because $F^\ttT$ is its left adjoint, $U^\ttT F^\ttT=T$, and the
corresponding comparison functor $\M^\ttT\to\M^\ttT$ sends the object $\bra
x,\alpha\ket$ to 
\[
\bra U^\ttT\bra x,\alpha\ket,U^\ttT\eps^\ttT_{\bra x,\alpha\ket}\ket
=\bra x,\alpha\ket\,.
\]
Therefore the comparison functor is the identity functor. Now we have an
equivalence $K\colon \C\to\M^\ttT$ and it is easy to see that
monadicity of $U^\ttT$ is inherited to $U$ via this $K$. Since $U=U^\ttT K$,
this defines a monoidal structure for $U$.
\end{proof}

For structural assumptions on $U$ and $\C$ that imply monadicity we refer to
the literature \cite{MacLane-Moerdijk, TTT}. Here we give only a crude
consequence of the above Theorem which is still sufficiently general to
include as special cases the forgetful functors of bialgebroids,
to be discussed in Section \ref{sec: bgd}. Therefore it covers also the cases
of forgetful functors $U\colon\,_A\M\to\,_kM$ where $A$ is either a weak
bialgebra or bialgebra or weak Hopf algebra or Hopf algebra over $k$.

\begin{cor}
Let $\C$ be a monoidal category having coequalizers and let $U\colon\C\to\M$
be a monoidal functor such that its underlying functor reflects
isomorphisms and has a left adjoint and a right adjoint. Then $\C$ is
monoidally equivalent to the Eilenberg-Moore category $\M^\ttT$ of a bimonad.
\end{cor}

Necessary and sufficient conditions for a bimonad to be the bimonad of
the forgetful functor of a bialgebroid will be given in Section \ref{sec: bgd}.

\section{2-functoriality of the construction of bimonads} \label{sec: 2func}

The functors $U$ for which a bimonad can be constructed are the objects of a
2-category $\LMF$. We extend the bimonad construction of the previous
Section to a 2-functor $\QQ\colon\LMF\to\Bmd$ from which a sensible
definition for the 2-category $\Bmd$ of bimonads emerges. We show that $\Q$ is
the left adjoint of a 2-functor $\EM$ which incorporates the Eilenberg-Moore
construction. This adjunction explains and extends the universality result
of Theorem \ref{thm: 1-univ}. Finally, the fact that bimonads
form a 2-category will enable us to speak about isomorphisms and equivalences
of bialgebroids which, in turn, in Section \ref{sec: bgd}, will be shown to
be objects of $\Bmd$.

\subsection{The 2-category of arrows}

Let $\bra \K,\b,\c\ket$ be a 2-category. As before in case of $\K=\Cat$ we omit
the symbol $\c$ for horizontal composition. We define the 2-category of arrows
in $\K$ as the 2-category $\Arr(\K)$ having
\begin{dinglist}{228}
\item objects $(\AA,U,\MM)$ where $U\colon \AA\to \MM$ is a 1-cell of $\K$,
\item 1-cells $\bra F,\kappa,G\ket\colon (\AA,U,\MM)\to(\BB,V,\NN)$ where
$F\colon \AA\to \BB$ and $G\colon \MM\to \NN$ are 1-cells of $\K$ and
$\kappa\colon GU\to VF$ is a 2-cell of $\K$,
\item 2-cells $[\vartheta,\nu]:\bra F,\kappa,G\ket\to\bra F',\kappa',G'\ket
\colon (\AA,U,\MM)\to(\BB,V,\NN)$ where \newline $\vartheta\colon F\to F'$ and
$\nu\colon G\to G'$ are 2-cells of $\K$ such that
\begin{equation}
V\vartheta\b\kappa=\kappa'\b \nu U
\end{equation}
\end{dinglist}
The horizontal composition of 1-cells is
\beq
\bra F,\kappa,G\ket\ohArr\bra H,\lambda,I\ket=\bra FH,\kappa H\b
G\lambda,GI\ket
\eeq
the horizontal composition of 2-cells is
\beq \label{eq: ohArr}
[\vartheta,\nu]\ohArr[\vartheta',\nu']=[\vartheta\vartheta',\nu\nu']
\eeq
and the vertical composition of 2-cells is
\beq \label{eq: ovArr}
[\vartheta,\nu]\ovArr[\vartheta',\nu']=[\vartheta\b\vartheta',\nu\b\nu']
\eeq
whenever they are defined.

Thus we have two 2-functors
\[\begin{CD}
\K@<\dom<<\Arr(\K)@>\cod>>\K
\end{CD}\]
given respectively by
\begin{eqnarray*}
\AA\mapsot&(\AA,U,\MM)&\mapsto \MM\\
F\mapsot&\bra F,\kappa,G\ket&\mapsto G\\
\vartheta\mapsot&[\vartheta,\nu]&\mapsto \nu
\end{eqnarray*}

We want to single out a sub-2-category in $\Arr(\K)$ the objects $(\AA,U,\MM)$
of which carry a universal action of a monoid $\bra T,\mu,\eta\ket$ at $\MM$.
Such actions $\alpha\colon TU\to U$ will be called left actions as they act on
the codomain side of $(\AA,U,\MM)$. As it is well known in universal algebra
\cite{Pareigis: NCAG} such monoids are readily obtained by universality from a
much simpler structure, a 1-cell $T\colon \MM\to \MM$ and a 2-cell
$\alpha\colon TU\to U$, no condition whatsoever. Existence of universal
monoids is guaranteed for example if $U$ has a left adjoint in $\K$. In the
next Definition the usual universality is replaced by a slightly stronger
"2-universality" property which we need later but which is also a property of
left adjointable $U$-s. In the sequel we denote by $\Arr^\sim(\K)$ the
sub-2-category in which the $\kappa$-s are invertible.

\begin{defi} \label{def: univ}
A left action on a 1-cell $(\AA,U,\MM)$ in $\K$ is a 2-cell in $\Arr(\K)$ of
the form $\bra \AA,\alpha,R\ket\colon (\AA,U,\MM)\to (\AA,U,\MM)$. That is to
say a left action on $U$ consists of a 1-cell $R\colon \MM\to \MM$ and a 2-cell
$\alpha\colon RU\to U$ in $\K$. The left action $\bfalpha=\bra
\AA,\alpha,R\ket$ is called universal if for any left action
$\bfbeta=\bra \BB,\beta, S\ket$ on $(\BB,V,\NN)$ and any 1-cells
$\bfkappa=\bra F,\kappa,G\ket$ and $\bfkappa'=\bra
F',\kappa',G'\ket$ from $(\AA,U,\MM)$ to $(\BB,V,\NN)$ in $\Arr^\sim(\K)$ the
domain functor gives rise to a bijection of 2-cells  
\beq
\begin{CD}
\Arr(\bfbeta\diamond\bfkappa,
\bfkappa'\diamond\bfalpha)@>\dom>>\K(F,F')\,. 
\end{CD}
\eeq
That is to say, a left action $\alpha\colon RU\to U$ is universal if for every 
left action $\beta\colon SV\to V$, 1-cells $F,F'\colon \AA\to \BB$, $G,G'\colon
\MM\to \NN$, 2-cell $\vartheta\colon F\to F'$ and invertible 2-cells $\kappa\colon
GU\iso VF$ and $\kappa'\colon G'U\iso VF'$ there exists a unique 2-cell
$\psi\colon SG\to G'R$ such that
\beq \label{eq: psi}
V\vartheta\b \beta F\b S\kappa = \kappa'\b G'\alpha\b \psi U\ .
\eeq
\end{defi}

\begin{rmk} \label{rmk: psi}
Taking into account the explicit form of the horizontal composites
\begin{align}
\bra \BB,\beta,S\ket\ohArr\bra F,\kappa,G\ket&=
\bra F,\beta F\b S\kappa,SG\ket\\
\bra F',\kappa',G'\ket\ohArr\bra \AA,\alpha,R\ket&=
\bra F',\kappa'\b G'\alpha,G'R\ket
\end{align}
the above equation (\ref{eq: psi}) for $\psi$ is precisely the condition for
the pair $[\vartheta,\psi]$ to be a 2-cell
\beq
[\vartheta,\psi]\colon\bra \BB,\beta,S\ket\ohArr\bra F,\kappa,G\ket
\to\bra F',\kappa',G'\ket\ohArr\bra \AA,\alpha,R\ket\,.
\eeq
Therefore, as sets
\beq
\{\psi\}=\cod_{\bfbeta\ohArr\bfkappa,
\bfkappa'\ohArr\bfalpha}\c
\dom^{-1}_{\bfbeta\ohArr\bfkappa,
\bfkappa'\ohArr\bfalpha}\ (\{\vartheta\})\,.
\eeq
\end{rmk}

The next Lemma secures a familiar class of
1-cells on which universal actions exist.

\begin{lem} \label{lem: L-adj=>univ}
If $U$ is a 1-cell in $\K$ which has a left adjoint $\bar U$ then there is a
universal left action on $U$, namely $\alpha=U\eps\colon TU\to U$ where
$T=U\bar U$ and $\eps\colon\bar U U\to \AA$ is the counit of the adjunction.
\end{lem}
\begin{proof}
Multiply equation (\ref{eq: psi}) from the left by ${\kappa'}^{-1}$, then
compose it horizontally from the right  with $\bar U$, and finally multiply it
from the right by $SG\eta$, where $\eta$ is the unit of the adjunction. Thus
we obtain
\begin{align*}
&{\kappa'}^{-1}\bar U\b V\vartheta\bar U\b\beta F\bar U\b S\kappa\bar U\b
SG\eta\\
=&G' U\eps\bar U\b \psi U\bar U\b SG\eta\\
=&G' U\eps\bar U\b G' U\bar U \eta\b\psi\\
=&\psi\ .
\end{align*}
\end{proof}

Let $\Arr_u(\K)$ denote the full sub-2-category of $\Arr(\K)$ having as objects
those objects of $\Arr(\K)$ on which a universal left action exists.

\subsection{The construction of the 2-functor $\Q$} \label{ss: Q}

We define a 2-functor $\Q\colon\Arr_u^\sim(\K)\to\Mnd(\K)$ as follows. 
Its object map is provided by the following Lemma.

\begin{lem} \label{lem: Q0}
For an object $(\AA,U,\MM)$ of $\Arr_u(\K)$ let $\bra \AA,\omega,T\ket$ be a
universal left action. Then there exist unique 2-cells $\mu\colon TT\to T$ and
$\eta\colon \MM\to T$ such that 
\begin{align} \label{eq: t-module1}
\omega\b\mu U&=\omega\b T\omega\\
\omega\b\eta U&=U\,. \label{eq: t-module2}
\end{align}
The triple $\bra T,\mu,\eta\ket$ is a monad in $\K$ on $\MM$.
\end{lem}
\begin{proof}
Let $\L(U)$ be the category of left actions on $U$ which is the subcategory in
$\Arr(U,U)$ containing as objects the special 1-cells $\bra \AA,\alpha,R\ket$
and as arrows the special 2-cells $[\AA,\nu]$. Then universality of $\bra
\AA,\omega,T\ket$ implies that it is a terminal object in $\L(U)$. As a matter
of fact if we specialize the universal property to the choice $V=U$, $\bra
F,\kappa,G\ket=\bra F',\kappa',G'\ket$ being the identity cell $(\AA,U,\MM)$
and $\vartheta=\AA$ then we obtain that for all $R\colon \MM\to \MM$ and all
$\alpha\colon RU\to U$ there exists a unique $\nu\colon R\to T$ such that
$\alpha=\omega\b \nu U$. 

Now it is clear that the solutions for $\mu$ and $\eta$ of the equations
(\ref{eq: t-module1}-\ref{eq: t-module2}) provide arrows 
\begin{align}
[\AA,\mu]\colon &\bra \AA,\omega\b T\omega,TT\ket\to \bra \AA,\omega,T\ket\\
[\AA,\eta]\colon&\bra \AA,U,\MM\ket\to\bra \AA,\omega,T\ket
\end{align}
hence they exist and are unique. The rest is standard universal algebra
\cite{Pareigis: NCAG}: One checks that both $\mu\b T\mu$ and $\mu\b \mu T$
provide arrows from $\bra \AA,\omega\b T\omega\b TT\omega,TTT\ket$ to the
terminal object, hence $\mu$ is associative. Similarly one proves that $\eta$
is a unit for $\mu$, hence $\bra T,\mu,\eta\ket$ is a monad.
\end{proof}

Given a choice of universal left action $\bra \AA,\omega,T\ket$ for each
object $(\AA,U,\MM)$ of $\Arr_u(\K)$ we define
\beq
\Q(\AA,U,\MM)\ :=\ \bra T,\mu,\eta\ket
\eeq
where the monad on the RHS is obtained from Lemma \ref{lem: Q0}.

The definition of $\Q$ on 1-cells is provided by the following specialization
of the universal property of Definition \ref{def: univ}. Setting $\bra
F,\kappa,G\ket=\bra F',\kappa',G'\ket$ and $\vartheta=F$ we obtain that if
$\bra \AA,\alpha,R\ket$ is a universal left action on $(\AA,U,\MM)$ then for
all object $(\BB,V,\NN)$, all left action $\bra \BB,\beta,S\ket$ on
$(\BB,V,\NN)$ and all $\bra F,\kappa,G\ket\colon (\AA,U,\MM)\to (\BB,V,\NN)$ in
$\Arr^\sim(\K)$ there exists a unique $\varphi\colon SG\to GR$ such that
\beq \label{eq: phi}
\beta F\b S\kappa\ =\ \kappa\b G\alpha\b\varphi U\,.
\eeq

\begin{lem} \label{lem: Q1}
Let $(\AA,U,\MM)$ and $(\BB,V,\NN)$ be objects of $\Arr_u(\K)$ and let
$\Q(\AA,U,\MM)$ $=\bra R,\mu_\ttR,\eta_\ttR\ket$ and $\Q(\BB,V,\NN)=\bra
S,\mu_\ttS,\eta_\ttS\ket$. Then for each 1-cell
\[
\bra F,\kappa,G\ket\colon (\AA,U,\MM)\to (\BB,V,\NN)\quad\in\ \Arr^\sim(\K)
\]
the unique
solution for $\varphi$ of equation (\ref{eq: phi}) provides a monad morphism
\beq
\bra G,\varphi\ket\colon\bra R,\mu_\ttR,\eta_\ttR\ket\to\bra
S,\mu_\ttS,\eta_\ttS\ket\,, 
\eeq
i.e., a 1-cell in the 2-category $\Mnd(\K)$ of monads in $\K$ \cite{Street:
Mnd}. That is to say $G\colon \MM\to \NN$ is a 1-cell and $\varphi\colon SG\to
GR$ is a 2-cell satisfying the commutative diagrams
\beq \label{dia: mnd mor}
\begin{picture}(160,70)
\put(20,50){$SSG$}
\put(43,53){\vector(1,0){35}} \put(53,56){$\scriptstyle S\varphi$}
\put(82,50){$SGR$}
\put(108,53){\vector(1,0){35}} \put(118,56){$\scriptstyle \varphi R$}
\put(148,50){$GRR$}
\put(30,45){\vector(0,-1){30}} \put(8,28){$\scriptstyle\mu_\ttS G$}
\put(22,5){$SG$}
\put(40,8){\vector(1,0){108}} \put(92,11){$\scriptstyle\varphi$}
\put(160,45){\vector(0,-1){30}} \put(164,28){$\scriptstyle G\mu_\ttR$}
\put(153,5){$GR$}
\end{picture}
\qquad\qquad
\begin{picture}(80,70)
\put(-5,5){$SG$}
\put(13,7){\vector(1,0){57}} \put(40,10){$\scriptstyle\varphi$}
\put(75,5){$GR$}
\put(39,50){$G$}
\put(35,45){\vector(-1,-1){30}} \put(5,35){$\scriptstyle\eta_\ttS G$}
\put(46,45){\vector(1,-1){30}} \put(62,35){$\scriptstyle G\eta_\ttR$}
\end{picture}
\eeq
\end{lem}
\begin{proof}
To prove the first diagram it suffices, by the above special universality 
property, to show that both $\varphi':=G\mu_\ttR\b\varphi R\b S\varphi$ and
$\varphi^\dagger:=\varphi\b\mu_\ttS G$ are solutions of
\beq
\beta^2F\b SS\kappa\ =\ \kappa\b G\alpha\b\Phi U
\eeq
for $\Phi\colon SSG\to GR$ where $\beta^2$ stands for $\beta\b
S\beta=\beta\b\mu_\ttS S\colon SSB\to B$. As a matter of fact
\begin{align*}
\kappa\b G\alpha\b\varphi' U
&=\kappa\b G\alpha\b G\mu_R U\b\varphi RU\b S\varphi U\\
&=\kappa\b G\alpha\b GR\alpha\b\varphi RU\b S\varphi U\\
&=\kappa\b G\alpha\b \varphi U\b SG\alpha\b S\varphi U\\
&=\beta F\b S\kappa\b SG\alpha\b S\varphi U\\
&=\beta F\b S\beta F\b SS\kappa\ =\ \beta^2F\b SS\kappa
\intertext{and}
\kappa\b G\alpha\b \varphi^\dagger U
&=\kappa\b G\alpha\b\varphi U\b\mu_S GU\\
&=\beta F\b S\kappa\b\mu_S GA\\
&=\beta F\b\mu_S VF\b SS\kappa\ =\ \beta^2F\b SS\kappa\,.
\end{align*}
Similarly, the proof of the second diagram amounts to showing that both
$G\eta_\ttR$ and $\varphi\b\eta_\ttS G$ solve the equation
\beq
\kappa\ =\ \kappa\b G\alpha\b\Upsilon U
\eeq
for $\Upsilon\colon G\to GR$. Indeed,
\begin{align*}
\kappa\b G\alpha\b G\eta_\ttR U&=\kappa\b GU\ =\ \kappa\\
\intertext{and}
\kappa\b G\alpha\b\varphi U\b\eta_\ttS GU&=\beta F\b S\kappa\b\eta_\ttS GU\\
&=\beta F\b\eta_\ttS VF\b\kappa\ =\ \kappa\,.
\end{align*}
\end{proof}

We can therefore define $\Q$ on 1-cells by
\beq
\Q\bra F,\kappa,G\ket\ :=\ \bra G,\varphi\ket
\eeq
where $\varphi$ is determined by Lemma \ref{lem: Q1}.

Before defining $\Q$ on 2-cells we investigate functoriality of $\Q$ on the
category of 0-cells and 1-cells. In the following Lemma $\alpha\colon RU\to U$,
$\beta\colon SV\to V$ and $\gamma\colon TW\to W$ denote universal actions
corresponding to the definition of $\Q$ on the object $U$, $V$ and $W$,
respectively.

\begin{lem} \label{lem: 1-functoriality of Q}
For composable 1-cells in $\Arr_u(\K)$ as in the diagram
\[
\begin{CD}
\AA@>F>>\BB@>H>>\CC\\
@VUV{\hskip 0.7truecm \kappa }V @V{V}V{\hskip 0.7truecm \lambda}V @VWVV\\
\MM@>G>>\NN@>I>>\OO
\end{CD}
\]
we have $\Q(\bra H,\lambda,I\ket\ohArr\bra F,\kappa,G\ket)=
\Q\bra H,\lambda,I\ket\ohBmd \Q\bra F,\kappa,G\ket$. 
\end{lem}
\begin{proof}
Taking into account the formula $\bra I,\chi\ket\ohBmd\bra G,\varphi\ket$=$\bra
IG,I\varphi\b\chi G\ket$ for composition of monad morphisms, we have to show
that if $\varphi\colon SG\to GR$ and $\chi\colon TI\to IS$ are solutions of
the equations
\begin{align*}
\beta F\b S\kappa&=\kappa\b G\alpha\b\varphi U\\
\gamma H\b T\lambda&=\lambda\b I\beta\b\chi V
\end{align*}
then $\nu=I\varphi\b\chi G$ solves the equation
\[
\gamma HF\b T\lambda F\b TI\kappa\ =\ \lambda F\b I\kappa\b IG\alpha\b\nu U\,.
\]
As a matter of fact,
\begin{align*}
&\lambda F\b I\kappa\b IG\alpha\b I\varphi U\b\chi GU=
\lambda F\b I\beta F\b IS\kappa\b\chi GU=\\
=&\lambda F\b I\beta F\b\chi VF\b TI\kappa=
\gamma HF\b T\lambda F\b TI\kappa\,.
\end{align*}
\end{proof}

The definition of $\Q$ on 2-cells uses the full strength of Definition
\ref{def: univ}. First of all for a 2-cell  
\[
[\vartheta,\nu]:\bra F,\kappa,G\ket\to\bra F',\kappa',G'\ket
\colon (\AA,U,\MM)\to(\BB,V,\NN)
\]
we set
\beq
\Q[\vartheta,\nu]\ :=\ \nu\,.
\eeq
The statement that $\nu\colon\bra G,\varphi\ket\to\bra G'\varphi'\ket$ is 
a transformation of monad morphisms, i.e., a 2-cell in $\Mnd(\K)$ is by
definition \cite{Street: Mnd} the property
\beq
\begin{CD}
SG@>\varphi>>GR\\
@V{S\nu}VV @VV{\nu R}V\\
SG'@>\varphi'>>G'R
\end{CD}
\eeq
where $\bra G,\varphi\ket=\Q\bra F,\kappa,G\ket$ and $\bra G',\varphi'\ket=\Q\bra
F',\kappa',G'\ket$. Commutativity of this diagram follows from universality
after noticing that both $\nu R\b\varphi$ and $\varphi'\b S\nu$ are solutions
for $\psi\colon SG\to G'R$ of the equation (\ref{eq: psi}). Indeed,
\begin{align*}
\kappa'\b G'\alpha\b \nu RU\b\varphi U&=
\kappa'\b\nu U\b G\alpha\b\varphi U\\
&=V\vartheta\b\kappa\b G\alpha\b\varphi U\ =\ V\vartheta\b\beta F\b S\kappa\\
\intertext{and}
\kappa'\b G'\alpha\b\varphi' U\b S\nu U&=
\beta F'\b S\kappa'\b S\nu U\\
&=\beta F'\b SV\vartheta\b S\kappa\ =\ V\vartheta\b\beta F\b S\kappa\,.
\end{align*}
Since for 2-cells both the horizontal and vertical composition in $\Mnd(\K)$
coincides with those of $\K$, in view of (\ref{eq: ohArr}) and (\ref{eq:
ovArr}) the $\Q$ preserves both compositions. This finishes the construction of
the 2-functor $\Q\colon\Arr^\sim_u(\K)\to\Mnd(\K)$.

\subsection{The monoidal version of $\Q$} \label{ss: QQ}

In this subsection we are interested in the monoidal properties of the
2-functor $\Q=\Q(\K)$ if the underlying 2-category $\K$ is monoidal. In order
not to drift too far from the main theme of bialgebroids we restrict ourselves
to the case $\K=\Cat$ endowed with the Cartesian product $\x$ of categories,
functors and natural transformations.

The content of this subsection crucially depends on whether $\alpha\colon
RU\to U$ being a universal action implies $\alpha\x\alpha\colon (R\x R)(U\x
U)\to U\x U$ is universal, too. Since this property does not seem to be
automatic, we shall restrict ourselves to functors $U$ with left adjoints. In
this case $U\x U$ also has a left adjoint, therefore $\alpha\x\alpha$ is
universal indeed (cf. Lemma \ref{lem: L-adj=>univ}). We denote by $\LFunc$
the full sub-2-category of $\Arr_u^\sim(\Cat)$ with objects the left
adjointable functors.

If $\A$ and $\M$ are monoidal categories and the left
adjointable $U\colon\A\to\M$ is given a monoidal structure $U_2\colon \oM(U\x
U)\iso U\oA$, $U_0\colon \iM\iso U\iA$ then Corollary \ref{cor: mon adj} 
tells us that the monad $\Q\bra\A,U,\M\ket=\bra T,\mu,\eta\ket$ has a lax
comonoidal structure $T^2\colon U\oA\to\oM(U\x U)$, $T^0\colon U\iA\to\iM$
so that $\bra T,T^2,T^0,\mu,\eta\ket$ is a bimonad. This result is the object
map part of a commutative diagram of 2-functors
\beq \label{dia: QQ}
\begin{CD}
\LMF@>\QQ>?>\Bmd\\
@VVV @VVV\\
\LFunc@>\Q>>\Mnd
\end{CD}
\eeq
where the vertical 2-functors forget about (co)monoidal structures, otherwise
all items in the first row are yet undefined, including the 2-category $\Bmd$
of bimonads. Our aim is to define them in such a way that the above diagram
commutes.

One solution is obtained by taking $\K$ to be $\ComonCat$ in the first row
and $\K=\Cat$ in the second and then apply the procedure of the previous
subsection to construct the $Q$. This seems to be the most natural choice
since bimonads involve lax comonoidal functors. This choice leads to the
diagram
\beq
\begin{CD}
\mathsf{L}\text{-}\ComonFunc@>\Q(\ComonCat)>>\Bmd'\\
@VVV @VVV\\
\LFunc@>\Q>>\Mnd
\end{CD}
\eeq
where $\mathsf{L}\text{-}\ComonFunc$ is the 2-category with
\begin{dinglist}{228}
\item objects the comonoidal functors $U$ (equivalently: monoidal ones) with
left adjoint, 
\item 1-cells $\bra F,\kappa, G\ket\colon U\to V$ where $F$, $G$ are lax
comonoidal functors and $\kappa\colon GU\to VF$ is a comonoidal natural
isomorphism,
\item 2-cells $[\vartheta,\nu]\colon \bra F,\kappa,G\ket\to\bra
F',\kappa',G'\ket$ where both $\vartheta\colon F\to F'$ and $\nu\colon G\to G'$
are comonoidal and satisfy the constraint $V\vartheta\b \kappa=\kappa'\b\nu U$.
\end{dinglist}

Accordingly, the 2-category $\Bmd'$ involves only (lax) comonoidal functors in
place of 1-cells. Especially, monad morphisms $\bra G,\varphi\ket$ involve lax
comonoidal functors $G\colon\M\to\N$. These functors map comonoids into
comonoids but does not map monoids to monoids. 

If we want arrows that preserve module algebras over bialgebroids instead of
module coalgebras, we must insist of having lax monoidal functors in the
definition of 1-cells. At first sight this spoils any sensible (co)monoidality
of the 2-cell $\varphi\colon SG\to GR$ since $R$ and $S$ are lax comonoidal
functors but $G$ is lax monoidal. Fortunately, the situation is not so bad.
\begin{defi} \label{def: ambi}
In the situation of the diagram
\[\begin{CD}
\M@>G>>\N\\
@VRVV @VVSV\\
\M'@>H>>\N'
\end{CD}\]
with four monoidal categories, lax monoidal functors $G$ and $H$ and lax
comonoidal functors $R$ and $S$ a natural transformation $\varphi\colon
SG\to HR$ is called ambimonoidal if the diagrams
\beq \label{dia: ambi1}
\begin{CD}
S\oN(G\x G)@>S^2(G\x G)>>\oNprime(SG\x SG)@>\oNprime(\varphi\x\varphi)>>
\oNprime(HR\x HR)\\
@V{SG_2}VV @. @VV{H_2(R\x R)}V\\
SG\oM@>\varphi\oM>>HR\oM@>HR^2>>H\oMprime(R\x R)
\end{CD}
\eeq
and
\beq \label{dia: ambi2}
\begin{picture}(160,70)
\put(20,50){$S\iN$}
\put(43,53){\vector(1,0){35}} \put(53,56){$\scriptstyle SG_0$}
\put(82,50){$SG\iM$}
\put(108,53){\vector(1,0){35}} \put(118,56){$\scriptstyle \varphi \iM$}
\put(148,50){$HR\iM$}
\put(30,45){\vector(0,-1){30}} \put(8,28){$\scriptstyle S^0$}
\put(22,5){$\iNprime$}
\put(40,8){\vector(1,0){108}} \put(92,11){$\scriptstyle H_0$}
\put(160,45){\vector(0,-1){30}} \put(164,28){$\scriptstyle HR^0$}
\put(153,5){$H\iMprime$}
\end{picture}
\eeq
are commutative.
\end{defi}

Beyond that it is meaningful the motivation for this definition comes from the
following
\begin{pro} \label{pro: ambi}
Let $U$ and $V$ be monoidal functors with left adjoints, $F$ and $G$ be lax
monoidal functors and let $\kappa\colon GU\to VF$ be a monoidal isomorphism
as in the diagram
\[
\begin{CD}
\A@>F>>\B\\
@VUV{\hskip 0.7truecm \kappa}V @VV{V}V\\
\M@>G>>\N
\end{CD}
\]
Forgetting the monoidal structures let $\Q(U)=R$, $\Q(V)=S$ and $\Q\bra
F,\kappa,G\ket=\bra G,\varphi\ket$. Then $\varphi\colon SG\to GR$ is
ambimonoidal if $R$ and $S$ are considered with lax comonoidal structures
by Theorem \ref{thm: F-sigma-pi}. 
\end{pro}
\begin{proof}
The defining equation (\ref{eq: phi}) of $\varphi$ is equivalent by Remark 
\ref{rmk: psi} to the condition that
\beq
[F,\varphi]\colon\bra\B,\beta,S\ket\ohArr\bra F,\kappa,G\ket\to
\bra F,\kappa,G\ket\ohArr\bra \A,\alpha,R\ket
\eeq
is a 2-cell in $\LFunc$. Monoidality of $\kappa$,
\beq
\kappa\oA\b GU_2\b G_2(U\x U)=VF_2\b V^2(F\x F)\b\oN(\kappa\x\kappa)
\eeq
is equivalent to
\beq
[F_2,G_2]\colon\bra\oB,V_2,\oN\ket\ohArr\left(\bra F,\kappa,G\ket\x\bra
F,\kappa,G\ket\right)\to
\bra F,\kappa,G\ket\ohArr\bra\oA,U_2,\oM\ket
\eeq
being a 2-cell in $\LFunc$. Comonoidality of $\beta=V\eps_V\colon SV\to V$,
\beq
V_2^{-1}\b\beta\oB=\oN(\beta\x\beta)\b S^2(V\x V)\b SV_2^{-1}
\eeq
is equivalent to
\beq
[\oB,S^2]\colon\bra\B,\beta,S\ket\ohArr\bra\oB,V_2,\oN\ket\to
\bra\oB,V_2,\oN\ket\ohArr\left(\bra\B,\beta,S\ket\x\bra\B,\beta,S\ket\right)
\eeq
being a 2-cell in $\LFunc$. Similarly, the $[\oA,R^2]$ is a 2-cell precisely
because $\alpha$ is comonoidal. Therefore one can take the following two
parallel vertical composites
\beq \label{eq: ambi1 LHS}
\begin{CD}
\bra\B,\beta,S\ket\ohArr\bra\oB,V_2,\oN\ket
\ohArr\left(\bra F,\kappa,G\ket\x\bra F,\kappa,G\ket\right)\\
@VV{\bra\B,\beta,S\ket\ohArr[F_2,G_2]}V\\
\bra\B,\beta,S\ket\ohArr\bra F,\kappa,G\ket\ohArr\bra\oA,U_2,\oM\ket\\
@VV{[F,\varphi]\ohArr\bra\oA,U_2,\oM\ket}V\\
\bra F,\kappa,G\ket\ohArr\bra\A,\alpha,R\ket\ohArr\bra\oA,U_2,\oM\ket\\
@VV{\bra F,\kappa,G\ket\ohArr[\oA,R^2]}V\\
\bra F,\kappa,G\ket\ohArr\bra\oA,U_2,\oM\ket\ohArr
\left(\bra\A,\alpha,R\ket\x\bra\A,\alpha,R\ket\right)
\end{CD}
\eeq
and
\beq \label{eq: ambi1 RHS}
\begin{CD}
\bra\B,\beta,S\ket\ohArr\bra\oB,V_2,\oN\ket
\ohArr\left(\bra F,\kappa,G\ket\x\bra F,\kappa,G\ket\right)\\
@VV{[\oB,S^2]\ohArr\left(\bra F,\kappa,G\ket\x\bra F,\kappa,G\ket\right)}V\\
\bra\oB,V_2,\oN\ket\ohArr\left(\bra\B,\beta,S\ket\x\bra\B,\beta,S\ket\right)
\ohArr\left(\bra F,\kappa,G\ket\x\bra F,\kappa,G\ket\right)\\
@VV{\bra\oB,V_2,\oN\ket\ohArr\left([F,\varphi]\x[F,\varphi]\right)}V\\
\bra\oB,V_2,\oN\ket\ohArr\left(\bra F,\kappa,G\ket\x\bra F,\kappa,G\ket\right)
\ohArr\left(\bra\A,\alpha,R\ket\x\bra\A,\alpha,R\ket\right)\\
@VV{[F_2,G_2]\ohArr\left(\bra\A,\alpha,R\ket\x\bra\A,\alpha,R\ket\right)}V\\
\bra F,\kappa,G\ket\ohArr\bra\oA,U_2,\oM\ket\ohArr
\left(\bra\A,\alpha,R\ket\x\bra\A,\alpha,R\ket\right)
\end{CD}
\eeq
Computing the $\dom$ of (\ref{eq: ambi1 LHS}) and (\ref{eq: ambi1 RHS}) we
obtain the same result, $F_2$. By universality of the left action
$\bfalpha\x\bfalpha$ their codomains should also be the same. Computing their
$\cod$ we obtain 
\bea
GR^2\b\varphi\oM\b SG_2\\
G_2(R\x R)\b \oN(\varphi\x\varphi)\b S^2(G\x G)
\eea
respectively. They are precisely the LHS and RHS of the ambimonoidality condition
(\ref{dia: ambi1}).

In order to prove the other ambimonoidality condition
\beq \label{eq: ambi2}
GR^0\b\varphi\iM\b SG_0\ =\ G_0\b S^0
\eeq
we note the following facts. The natural transformations
\beq
\begin{CD}
1@>\iA>>\A\\
@| @V{U_0\hskip 0.5 truecm}V{U}V\\
1@>\iM>>\M
\end{CD}
\qquad\text{and}\qquad
\begin{CD}
1@>\iB>>\B\\
@| @V{V_0\hskip 0.5 truecm}V{V}V\\
1@>\iN>>\N
\end{CD}
\eeq
where $1$ is the one element category, are 1-cells in $\LFunc$ because $U_0$
and $V_0$ are invertible. Counitality of the comonoidal natural transformation
$\alpha$,
\beq
U_0^{-1}\b\alpha\iA\ =\ R^0\b RU_0^{-1}
\eeq
is equivalent to the statement that
\beq
[\iA,R^0]\colon\bra\A,\alpha,R\ket\ohArr\bra\iA,U_0,\iM\ket\to
\bra\iA,U_0,\iM\ket\colon(1,1,1)\to(\A,U,\M)
\eeq
is a 2-cell in $\LFunc$. Similarly, counitality of $\beta$ is the condition
for $[\iB,S^0]$ to be a 2-cell.
Unitality of the monoidal natural transformation $\kappa$,
\beq
VF_0\b V_0\ =\ \kappa\iA\b GU_0\b G_0
\eeq
in turn is the condition for
\beq
[F_0,G_0]\colon\bra\iB,V_0,\iN\ket\to 
\bra F,\kappa,G\ket\ohArr\bra\iA,U_0,\iM\ket
\eeq
to be a 2-cell in $\LFunc$. Thus one may form the vertical composites
\beq
\begin{CD}
\bra\B,\beta,S\ket\ohArr\bra\iB,V_0,\iN\ket\\
@VV{ \bra\B,\beta,S\ket\ohArr[F_0,G_0]}V\\
\bra\B,\beta S\ket\ohArr \bra F,\kappa,G\ket\ohArr\bra\iA,U_0,\iM\ket\\
@VV{[F,\varphi]\ohArr\bra\iA,U_0,\iM\ket}V\\
\bra F,\kappa,G\ket\ohArr\bra\A,\alpha,R\ket\ohArr\bra\iA,U_0,\iM\ket\\
@VV{\bra F,\kappa,G\ket\ohArr[\iA,R^0]}V\\
\bra F,\kappa,G\ket\ohArr\bra\iA,U_0,\iM\ket
\end{CD}
\eeq
and
\beq
\begin{CD}
\bra\B,\beta,S\ket\ohArr\bra\iB,V_0,\iN\ket\\
@VV{[\iB,S^0]}V\\
\bra \iB,V_0,\iN\ket\\
@VV{[F_0,G_0]}V\\
\bra F,\kappa,G\ket\ohArr\bra\iA,U_0,\iM\ket
\end{CD}
\eeq
which evaluate to be
\bea
&[F_0,GR^0\b\varphi\iM\b SG_0]\\
&[F_0,G_0\b S^0]
\eea
respectively. Since their $\dom$ are equal, we can conclude by universality
that their $\cod$ are equal as well, which in turn are the LHS and RHS of
(\ref{eq: ambi2}). The universal action we use here is the trivial action
$(1,1,1)$ on the identity functor of the category $1$. The cell $U_0$, using
its invertibility, can be absorbed into $\bfkappa\ohArr U_0$ to form the
$\bfkappa'$ of the universality condition of Definition \ref{def: univ}.
Universality of $(1,1,1)$ in turn follows directly from (\ref{eq: psi})
noticing that after inserting a 1-cell for $\alpha$ and a 0-cell for $U$
equation (\ref{eq: psi}) immediately gives a unique solution for $\psi$.
\end{proof}

Motivated by the above Proposition we can now fill in the missing items in
diagram (\ref{dia: QQ}). 
\begin{defi} \label{def: MonFunc_0}
Let $\LMF$ be the 2-category with
\begin{dinglist}{228}
\item objects the monoidal functors $U\colon\A\to\M$ with the underlying
functor having a left adjoint,
\item 1-cells $(\A,U,\M)\to (\B,V,\N)$ the triples $\bra F,\kappa,G\ket$
where $F\colon\A\to\B$, $G\colon\M\to\N$ are lax monoidal functors and
$\kappa\colon GU\to VF$ is a monoidal natural isomorphism and
\item 2-cells $\bra F,\kappa,G\ket\to\bra F',\kappa',G'\ket$ the pairs
$[\vartheta,\nu]$  where $\vartheta \colon F\to F'$ and $\nu\colon G\to G'$ are
monoidal natural transformations satisfying the
constraint $V\vartheta\b\kappa=\kappa'\b\nu U$. 
\end{dinglist}
All compositions are defined as in $\LFunc$ via forgetting $\LMF\to
\LFunc$.
\end{defi}
This 2-category describes the precise framework in which we are able to
associate a "quantum groupoid" to a forgetful functor. The "quantum groupoids"
in this generality are the bimonads.
\begin{defi}
Let $\Bmd$ be the 2-category with
\begin{dinglist}{228}
\item objects the bimonads $\bra \M,\ttT\ket\equiv\bra \M,T,T^2,T^0,
\mu,\eta\ket$ of Definition \ref{def: bmd},
\item 1-cells $\bra\M,\ttR\ket\to \bra\N,\ttS\ket$ the monad morphisms $\bra
G,\varphi\ket$ in which the functor $G\colon \M\to \N$ is lax monoidal 
and $\varphi\colon SG\to GR$ is ambimonoidal,
\item 2-cells $\bra G,\varphi\ket\to\bra G',\varphi'\ket$ the monad
transformations $\nu\colon G\to G'$ which are monoidal.
\end{dinglist}
All compositions are defined by the forgetting 2-functor $\Bmd\to\Mnd$.
\end{defi}
In order for $\Bmd$ to be well defined as a 2-category we are still indebted to
show that ambimonoidality is preserved by horizontal composition.
\begin{lem}
Consider two horizontally composable monad morphisms
\[
\begin{CD}
\bra\A,R\ket@>\bra G,\varphi\ket>>\bra\B,S\ket@>\bra H,\chi\ket>>\bra \C,T\ket
\end{CD}
\]
in which the categories are monoidal, the functors are lax monoidal and the
natural transformations $\varphi$ and $\chi$ are ambimonoidal. Then in the
composite
\[
\bra H,\chi\ket\ohBmd\bra G,\varphi\ket\ =\ \bra HG,H\varphi\b\chi G\ket
\]
the functor is lax monoidal and the natural transformation is ambimonoidal.
\end{lem}
\begin{proof}
Lax monoidality of $HG$ is obvious. We have to show that 
$\xi:=H\varphi\b\chi G$ satisfies the two diagrams (\ref{dia: ambi1}) and
(\ref{dia: ambi2}). The proof is this. The first ambimonoidality diagram
follows from commutativity of

\beq
\begin{CD}
T\o(HG\x HG)@>TH_2(G\x G)>>TH\o(G\x G)@>THG_2>>THG\o\\
@VV{T^2(HG\x HG)}V @VV{\chi\o(G\x G)}V @VV{\chi G\o}V\\
\o(THG\x THG) @. HS\o(G\x G)@>HSG_2>>HSG\o\\
@VV{\o(\chi G\x \chi G)}V @VV{HS^2(G\x G)}V @VV{H\varphi\o}V\\
\o(HSG\x HSG)@>>H_2(SG\x SG)>H\o(SG\x SG) @. HGR\o\\
@VV{\o(H\varphi\x H\varphi)}V @VV{H\o(\varphi\x\varphi)}V @VV{HGR^2}V\\
\o(HGR\x HGR)@>>H_2(GR\x GR)>H\o(GR\x GR)@>>HG_2(R\x R)>HG\o(R\x R)
\end{CD}
\eeq
where $\o$ in the 1st, 2nd and 3rd column denotes the monoidal product of
$\C$, $\B$ and $\A$, respectively. The second ambimonoidality diagram follows
from

\beq
\begin{picture}(200,120)
\put(0,10){$\iC$}  \put(60,10){$H\iB$} \put(210,10){$HG\iA$}
\put(60,60){$HS\iB$} \put(130,60){$HSG\iA$} \put(210,60){$HGR\iA$}
\put(0,110){$T\iC$} \put(60,110){$TH\iB$} \put(130,110){$THG\iA$}
\put(10,12){\vector(1,0){45}}  \put(30,15){$\scriptstyle H_0$}
\put(85,12){\vector(1,0){120}} \put(130,15){$\scriptstyle HG_0$}
\put(90,62){\vector(1,0){35}}  \put(92,65){$\scriptstyle HSG_0$}
\put(170,62){\vector(1,0){35}} \put(177,65){$\scriptstyle H\varphi\iA$}
\put(15,112){\vector(1,0){40}} \put(30,115){$\scriptstyle TH_0$}
\put(90,112){\vector(1,0){35}} \put(92,115){$\scriptstyle THG_0$}
\put(5,105){\vector(0,-1){85}} \put(-10,60){$\scriptstyle T^0$}
\put(68,105){\vector(0,-1){35}} \put(50,85){$\scriptstyle \chi\iB$}
\put(68,55){\vector(0,-1){35}}  \put(48,35){$\scriptstyle HS^0$}
\put(140,105){\vector(0,-1){35}} \put(145,85){$\scriptstyle \chi G\iA$}
\put(220,55){\vector(0,-1){35}} \put(225,35){$\scriptstyle HGT^0$}
\end{picture}
\eeq
\end{proof}

Let us summarize what we have obtained sofar:
\begin{thm}
Given a 2-functor $\Q\colon \LFunc\to\Mnd$ as in Subsection \ref{ss: Q}
there is a unique 2-functor
\[
\QQ\colon\LMF\to\Bmd
\]
such that (\ref{dia: QQ}) is commutative. 
\end{thm}
\begin{proof}
As for the object map of $\QQ$ we must take the monoid $\Q(U)$ and endow it
with the comonoidal structure that Corollary \ref{cor: mon adj} provides.
The arrow map of $\QQ$ is again uniquely determined by that of $\Q$ and it
yields bimonad morphisms by Proposition \ref{pro: ambi}. For the unique 2-cell
map of $\QQ$ there is nothing to prove.
\end{proof}

\subsection{The monoidal Eilenberg-Moore construction as a 2-functor}

Functoriality of the Eilenberg-Moore construction can be formalized as having
a 2-functor the object map of which associates forgetful functors to bimonads.

Let $\EM\colon\Bmd\to\LMF$ be the 2-functor defined as follows.

\textit{The object map:} For a bimonad $\bra \M,\ttT\ket$ let
$\EM\bra\M,\ttT\ket:=(\M^\ttT,U^\ttT,\M)$, the strict monoidal forgetful
functor of the category $\M^\ttT$ of $\ttT$-algebras, see Proposition
\ref{pro: mon EM}. 

\textit{The arrow map:} For a bimonad morphism $\bra G,\varphi\ket$ we
define 
\[
\EM\left(\bra\M,\ttT\ket\rarr{\bra G,\varphi\ket}\bra
\N,\ttS\ket\right)\ =\ \bra G^\varphi,=,G\ket
\]
where 
\[
\bra G^\varphi,=,G\ket\ =\ 
\begin{CD}
\M^\ttT@>G^\varphi>>\N^\ttS\\
@VV{U^\ttT\hskip 0.5truecm =}V @VV{U^\ttS}V\\
\M@>G>>\N
\end{CD}
\qquad\quad\text{and}\quad
G^\varphi\colon
\begin{CD}
\bra x,\alpha\ket@.\ \mapsto\ @.\bra Gx,G\alpha\c\varphi_x\ket\\
@VV{\tau \hskip 0.5truecm\mapsto}V @. @VV{G\tau}V\\
\bra y,\beta\ket@.\ \mapsto\ @.\bra Gy,G\beta\c\varphi_y\ket
\end{CD}
\]
which is indeed a functor from $\ttT$-algebras to $\ttS$-algebras since
$G\beta\c\varphi_y\c SG\tau=G\beta\c GT\tau\c\varphi_x=G\tau\c
G\alpha\c\varphi_x$. The monoidal structure for $G^\varphi$ is the one
given in Lemma \ref{lem: mon G^fi} below.

\textit{The 2-cell map:} For a transformation $\nu\colon \bra
G,\varphi\ket\to \bra G',\varphi'\ket$ of monad morphisms $\bra
\M,\ttT\ket\to\bra\N,\ttS\ket$ we define
\[
\EM(\nu)\ :=\ [\hat\nu,\nu]\colon\bra G^\varphi,=,G\ket\to\bra
G'^{\varphi'},=,G'\ket
\]
where $\hat\nu$ on the $\ttT$-algebra $\bra x,\alpha\ket$ is the
lift of $\nu_x$,
\[
\hat\nu_{\bra x,\alpha\ket}=\left(
\begin{CD}
\bra Gx,G\alpha\c\varphi_x\ket@>\nu_x>>\bra G'x,G'\alpha\c\varphi'_x\ket
\end{CD}
\right)
\]
which is indeed an $\ttS$-algebra morphism because
$G'\alpha\c\varphi'_x\c S\nu_x$ $=G'\alpha\c\nu_{Tx}\c\varphi_x$ $=\nu_x\c
G\alpha\c\varphi_x$.

\begin{lem} \label{lem: mon G^fi}
\begin{align*}
G^\varphi_2&=\left\{
\begin{CD}
G^\varphi\bra x,\alpha\ket\oalg_\N G^\varphi\bra y,\beta\ket
@>G_{x,y}>>
G^\varphi(\bra x,\alpha\ket\oalg_\M \bra y,\beta\ket)
\end{CD}
\right\}\\
G^\varphi_0&=\left(
\begin{CD}
\bra\iN,S^0\ket
@>G_0>>
G^\varphi\bra\iM,T^0\ket
\end{CD}
\right)
\end{align*}
is the unique monoidal structure on $G^\varphi$ such that $U^\ttS G^\varphi=
G U^\ttT$, as monoidal functors.
\end{lem}
\begin{proof}
Since $U^\ttT$ and $U^\ttS$ are strict monoidal, the only monoidal structure
on $G^\varphi$ is the one with components that are lifted from the components
of $G_2$, $G_0$, which is precisely the above formula. The hexagon and square
identities therefore hold automatically if we can show that the components
$G_{x,y}$ and $G_0$ can indeed be lifted to $\ttS$-algebra maps.

$G_{x,y}$ lifts to an arrow in $\N^\ttS$ iff 
\[
G(\alpha\oM\beta)\c GT_{x,y}\c\varphi_{x\oM y}\c SG_{x,y}
=G_{x,y}\c(G\alpha\o_N G\beta)\c(\varphi_x\oN\varphi_y)\c S_{Gx,Gy}
\]
which, after using naturality of $G_2$, becomes a consequence of the first
ambimonoidality axiom for $\varphi$.

$G^\varphi_0$ is an arrow in $\N^\ttS$ iff $GT^0\c\varphi_{\iM}\c SG_0=G_0\c
T^0$ which is precisely the second ambimonoidality axiom for $\varphi$.
\end{proof}

\subsection{The adjunction $\QQ\dashv\EM$ and universality}

In this subsection we construct pseudo natural transformations $\xi$ and
$\zeta$ providing the unit and counit of the adjunction $\QQ\dashv\EM$,
respectively. Then we show how to restrict these 2-functors to obtain
an adjunction in the strict sense, which is needed to establish the universal
property of the bimonad $\QQ(U)$ of a left adjointable monoidal functor. 

\subsubsection{The counit $\zeta$}
The action of $\QQ$ on $\EM(\M,T)$ does not necessarily return the original
monad $T$. Instead it gives $(\M,U^\ttT\bar U^\ttT)$ where $\bar U^\ttT$ is
some left adjoint of $U^\ttT$. If $\bar U^\ttT$ were equal to the free
$\ttT$-algebra functor $F^\ttT$ then we would get the original monad $\ttT$.
Of course, all left adjoints are isomorphic and it is easy to see that the
isomorphism $\sigma\colon F^\ttT\iso\bar U^\ttT$ leads to a monad morphism
$\zeta=U^\ttT\sigma$ sending $T=U^\ttT F^\ttT$ to $U^\ttT\bar U^\ttT$. A
closer look gives that we are actually have a bimonad isomorphism. This is the
content of the next Lemma in which $\twocom$ denotes composition of
2-functors.  
\begin{lem} \label{lem: zeta1}
$\bra\M,\zeta\ket\colon\QQ\twocom\EM(\M,T)\iso(\M,T)$ is a bimonad isomorphism.
\end{lem}
\begin{proof}
Because of uniqueness of left adjoints of monoidal functors up to
comonoidal natural isomorphisms (Lemma \ref{lem: uniq mon adj}), the $\sigma$
can be chosen to be comonoidal. But $U^\ttT$ is also (co)monoidal, so the
$\zeta$ is, either. Now the identity functor $\M$ being comonoidal, the
ambimonoidality condition for $\zeta$ is equivalent to its comonoidality. 
\end{proof}
Next we investigate the naturality properties of $\zeta$. 
Let $\bra G,\varphi\ket$ be a bimonad morphism $(\M,\ttT)\to(\N,\ttS)$. Then
$\QQ\twocom\EM\bra G,\varphi\ket=\bra G,\varphi'\ket$ where $\varphi'$ is
determined from the 1-cell 
\[
\EM\bra G,\varphi\ket\ =\ 
\begin{CD}
\M^\ttT@>G^\varphi>>\N^\ttS\\
@VV{U^\ttT\hskip 0.5truecm}V @VV{U^\ttS}V\\
\M@>G>>\N
\end{CD}
\]
i.e., $\varphi\colon S'G\to GT'$ is the unique solution of
\[
\beta' G^\varphi\ =\ G\alpha'\b\varphi' U^\ttT
\]
where $T'=U^\ttT\bar U^\ttT$, $S'=U^\ttS\bar U^\ttS$ and
$\alpha'=U^\ttT\eps_{U^\ttT}$, $\beta'= U^\ttS\eps_{U^\ttS}$ are the universal
actions associated to $U^\ttT$ and $U^\ttS$, respectively, in the
definition of $\QQ$. Denoting by $\alpha=U^\ttT\eps^\ttT$ and
$\beta=U^\ttS\eps^\ttS$, respectively, the universal actions associated to
them by the Eilenberg-Moore construction, we have 
\[
\alpha'=\alpha\b\zeta_\ttT^{-1}U^\ttT\,,\qquad
\beta'=\beta\b\zeta_\ttS^{-1}U^\ttS
\]
so we have to solve
\[
\beta G^\varphi\b\zeta_\ttS^{-1}U^\ttS G^\varphi\ =\ G\alpha\b
G\zeta_\ttT^{-1} U^\ttT\b\varphi'U^\ttT\,.
\]
Taking into account the formulae below in which $\bra x,\alpha\ket$ stands for
any $\ttT$-algebra 
\begin{align*}
\eps^\ttS_{G^\varphi\bra x,\alpha\ket}&=
\eps^\ttS_{\bra Gx,G\alpha\c\varphi_x\ket}=\bra
SGx,\mu^\ttS_{Gx}\ket\rarr{G\alpha\c\varphi_x}\bra Gx,G\alpha\c\varphi_x\ket\\
\beta_{G^\varphi\bra x,\alpha\ket}&=SGx\rarr{G\alpha\c\varphi_x} Gx\\
G\alpha_{\bra x,\alpha\ket}&= GTx\rarr{G\alpha} Gx
\end{align*}
the solution is
\[
\varphi'=G\zeta_\ttT\b\varphi\b\zeta_\ttS^{-1}G
\]
which is equivalent to that
\beq
\begin{CD}
(\M,\ttT')@>\bra\M,\zeta_\ttT\ket>>(\M,\ttT)\\
@V{\bra G,\varphi'\ket}V{\hskip 0.7truecm =}V @VV{\bra G,\varphi\ket}V\\
(\N,\ttS')@>\bra\N,\zeta_\ttS\ket>>(\N,\ttS)
\end{CD}
\eeq
is an identity 2-cell.
\begin{lem} \label{lem: zeta2}
$\zeta\colon\QQ\twocom\EM\iso\Bmd$ is a 2-natural isomorphism, i.e., for any
2-cell $\nu\colon\bra G,\varphi\ket\to \bra G',\varphi'\ket\colon
(\M,\ttT)\to(\N,\ttS)$ 
\[
\nu\ohBmd\zeta_\ttT\ =\ \zeta_\ttS\ohBmd\QQ\twocom\EM(\nu)\,.
\]
\end{lem}
\begin{proof}
In the above preparations we have already shown this relation for 1-cells
$\nu$. If $\nu$ is a 2-cell then it suffices to check the equation merely as
an equality of natural transformations, i.e., as 2-cells in $\Cat$. Since the
functor component of $\zeta$ is always the identity functor, this equality is
the trivial $\nu=\nu$.
\end{proof}

\subsubsection{The unit $\xi$}

The Eilenberg-Moore comparison functors $K_U\colon \A\to\M^\ttT$ provide
1-cells
\beq
\xi_U:=\bra K_U,=,\M\ket\ =\ 
\begin{CD}
\A@>K_U>>\M^\ttT\\
@V{U}V{\hskip 0.7truecm}V @VV{U^\ttT}V\\
\M@=\M
\end{CD}
\quad :\ U\to\EM\twocom\QQ(U)
\eeq
for all objects $(\A,U,\M)$ in $\LMF$. On 1-cells the 2-functor 
$\EM\twocom\QQ$ acts as
\beq
\begin{CD}
(\A,U,\M)\\
@VV{\bra F,\kappa,G\ket}V\\
(\B,V,\N)
\end{CD}
\qquad\qquad
\stackrel{\EM\twocom\QQ}{\mapsto}\quad
\begin{CD}
(\M^\ttT,U^\ttT,\M)\\
@VV{\bra G^\varphi,=,G\ket}V\\
(\N^\ttS,U^\ttS,\N)
\end{CD}
\eeq
where $\varphi\colon SG\to GT$ is the unique solution of
\beq
\beta F\b S\kappa\ =\ \kappa\b G\alpha\b\varphi U
\eeq
where $T=U\bar U$, $S=V\bar V$, $\alpha=U\eps_U$ and $\beta=V\eps_V$.
\begin{lem} \label{lem: xi1}
The monoidal natural isomorphism $\kappa\colon GU\to VF$ lifts to a monoidal
natural isomorphism $\hat\kappa\colon G^\varphi K_U\to K_VF$. The pair
$\Xi_\kappa:=[\hat\kappa,\id]$ is an invertible 2-cell
\beq
\begin{CD}
U@>\xi_U>>\EM\twocom\QQ(U)\\
@V{\kappa}V{\hskip 0.5truecm\Xi_\kappa \swarrow}V
@VV{\EM\twocom\QQ(\kappa)}V\\
V@>>\xi_V>\EM\twocom\QQ(V)
\end{CD}
\eeq
in $\LMF$.
\end{lem}
\begin{proof}
Computing the effect of the functors on an object $a\in\A$
\begin{align*}
K_VF\colon &a\mapsto\bra VFa,SVFa\rarr{\beta Fa}VFa\ket\\
G^\varphi K_U\colon &a\mapsto\bra GUa,SGUa\rarr{\varphi Ua}GRUa\rarr{G\alpha
a}GUa\ket
\end{align*}
we see that the lifting property of $\kappa_a\colon GUa\to VFa$ is just the
defining equation of $\varphi$ above. So $\hat\kappa$ has the proper
components and it is natural by virtue of the very simple form of the functors
$K_U$, $K_V$ on arrows. Moreover, $\hat\kappa$ is monoidal since the monoidal
structures of $G^\varphi$, $K_U$, $K_V$ are just the lifts of the
corresponding structures in $G$, $U$, $V$, respectively. The constraint for
$\Xi_\kappa=[\hat\kappa,\id]$ to be a 2-cell is just the lifting property
$U^\ttS\hat\kappa=\kappa$.
\end{proof}
  
\begin{lem} \label{lem: xi2}
$\xi\colon\LMF\to\EM\twocom\QQ$ is a pseudo natural transformation with 
\[
\Xi_\kappa\colon \EM\twocom\QQ(\kappa)\ohArr\xi_U\iso\xi_V\ohArr\kappa\,,
\qquad \kappa\colon U\to V\,.
\]
That is to say, for any 2-cell $\theta\colon\kappa\to\kappa'\colon U\to V$ 
\[
\Xi_{\kappa'}\ovArr(\EM\twocom\QQ(\theta)\ohArr\xi_U)
=(\xi_V\ohArr\theta)\ovArr\Xi_\kappa\ .
\]
\end{lem}
\begin{proof}
For $\theta=[\vartheta,\nu]$ one has $\EM\twocom\QQ(\theta)=[\hat\nu,\nu]$ so
one has to check $\dom$ and $\cod$ of
\[
[\hat\kappa,=]\ovArr\left([\hat\nu,\nu]\ohArr\bra K_U,=,\M\ket\right)=
\left(\bra K_V,=,\N\ket\ohArr[\vartheta,\nu]\right)\ovArr[\hat\kappa,=]
\]
which are $\hat\kappa'\b\hat\nu K_U=K_V\vartheta\b\hat\kappa$ and the identity
$\nu=\nu$, respectively. The former is the lift of $\kappa'\b\nu
U=V\vartheta\b\kappa$ which is but the the defining equation for the pair
$[\vartheta,\nu]$ to be a 2-cell $\kappa\to\kappa'$.
\end{proof}

\subsubsection{The pseudo adjunction $\QQ\dashv\EM$}
We want to prove that the 2-functor $\QQ$ is the pseudo left adjoint of $\EM$
in the following sense.
\begin{thm}
There exist pseudo natural transformations 
\[
\zeta\colon \QQ\twocom\EM\to\Bmd\,,\qquad\xi\colon\LMF\to\EM\twocom\QQ
\]
such that
\begin{align}
(\zeta\twocom\QQ)\onecom(\QQ\twocom\xi)&=\QQ\label{eq: psadj1}\\
(\EM\twocom\zeta)\onecom(\xi\twocom\EM)&=\EM\label{eq: psadj2}
\end{align}
where $\twocom$ denotes composition of 2-functors and higher cells, i.e.,
2-composition, and $\onecom$ denotes componentwise horizontal composition of
natural transformations, i.e., 1-composition in the 3-category 2-$\Cat$.
\end{thm}
\begin{proof}
We use the pseudo natural transformations $\zeta$ and $\xi$ constructed in 
the Lemmas \ref{lem: zeta1}, \ref{lem: zeta2}, \ref{lem: xi1} and \ref{lem:
xi2}. The effect of the 2-functor $\QQ$ on the 1-cell $\xi_U$ is
\[
\QQ(\xi_U)=\QQ\bra
K_U,=\M\ket=\left((\M,T)\rarr{\bra\M,\zeta'\ket}(\M,T')\right)
\]
where $T=U\bar U$, $T'=U^\ttT\bar U^\ttT$ and $\zeta'\colon T'\M\to\M T$ is
the unique solution of
\[
\alpha' K_U=\alpha\b\zeta'U
\]
where $\alpha=U\eps_U$, and
$\alpha'=U^\ttT\eps_{U^\ttT}=U^\ttT\eps^\ttT\b\zeta_\ttT^{-1}U^\ttT$.
Since
\begin{align*}
\alpha'K_Ua&=\alpha'_{\bra Ua,U\eps_Ua\ket}=U\eps_Ua\c\zeta_\ttT^{-1}Ua\\
\alpha'K_U&=\alpha\b\zeta_\ttT^{-1}U\,,
\end{align*}
the solution is $\zeta'=\zeta_\ttT^{-1}$. Since $\ttT$ here means $\QQ(U)$, we
obtain
\[
\left(\QQ(U)\rarr{\QQ(\xi_U)}\QQ\twocom\EM\twocom\QQ(U)
\rarr{\zeta_{\QQ(U)}}\QQ(U)\right)\ =\ 1_{\QQ(U)}
\]
which is precisely equation (\ref{eq: psadj1}). In order to prove the other
adjunction relation we look at
\[
\xi_{\EM(T)}=\bra K_{U^\ttT},=,\M\ket\colon (\M^\ttT,U^\ttT,\M)\to
(\M^{\ttT'},U^{\ttT'},\M)
\]
where $T'=U^\ttT\bar U^\ttT$ as before and
\[
K_{U^\ttT}\colon\bra x,\alpha\ket\mapsto\bra U^\ttT\bra x,\alpha\ket,
U^\ttT\eps_{U^\ttT}\bra x,\alpha\ket\ket=\bra x,\alpha\c\zeta_x^{-1}\ket
\ .\]
Let us compare this with
\[
\EM\left((\M,T')\rarr{\zeta}(\M,T)\right)\ =\ 
\begin{CD}
\M^{\ttT'}@>\M^\zeta>>\M^\ttT\\
@VV{U^{\ttT'}\hskip 0.5truecm =}V @VV{U^\ttT}V\\
\M@=\M
\end{CD}
\]
where $\M^\zeta\colon\bra x,\alpha\ket\mapsto\bra x,\alpha\c\zeta_x\ket$, i.e.,
$\M^\zeta=(K_{U^\ttT})^{-1}$. Therefore the 1-cell $\EM(\zeta_\ttT)$ is the
strict inverse of $\xi_{\EM(T)}$ for any bimonad $T$. This proves 
\[
\left(\EM(T)\rarr{\xi_{\EM(T)}}\EM\twocom\QQ\twocom\EM(T)\rarr{\EM(\zeta_\ttT)}
\EM(T)\right)\ =\ 1_{\EM(T)}
\] 
for all bimonad $T$ which is equation (\ref{eq: psadj2}).
\end{proof}

\subsubsection{Universality}

The fact that $\QQ$ is a left pseudo adjoint of $\EM$ has the following local
description. For each object $U$ in $\LMF$ there exist a bimonad
$T=\QQ(U)$ and a 1-cell $\xi=\xi_U\colon U\to\EM(T)$ satisfying the following
property:
\begin{description}
\item[P] If $S$ is a bimonad and $\kappa\colon U\to\EM(S)$ is a 1-cell then
there exists a, up to isomorphism unique, bimonad morphism $\varphi\colon T\to
S$ such that 
\[
\EM(\varphi)\ohArr\xi\cong\kappa\,.
\]
\end{description}
As a matter of fact, let $\varphi:=\zeta_S\ohBmd\QQ(\kappa)$. Then
\begin{align*}
\EM(\varphi)\ohArr\xi&=\EM(\zeta_S)\ohArr\EM\twocom\QQ(\kappa)\ohArr\xi\iso\\
&\begin{CD}@>\EM(\zeta_S)\ohArr\Xi_\kappa>\sim>
\EM(\zeta_S)\ohArr\xi_{\EM(S)}\ohArr\kappa=\kappa\,.
\end{CD}
\end{align*}
If $\varphi'\colon T\to S$ is another monad morphism for which there exists a
$\phi'\colon\EM(\varphi')\ohArr\xi\iso\kappa$ then
\begin{align*}
\varphi'&=\varphi'\ohBmd\zeta_{\QQ(U)}\ohBmd\QQ(\xi)\rarr{=}\zeta_S\ohBmd
\QQ\twocom\EM(\varphi')\ohBmd\QQ(\xi)\iso\\
&\begin{CD}@>\zeta_S\ohBmd\QQ(\phi')>\sim>\zeta_S\ohBmd\QQ(\kappa)=\varphi
\end{CD}
\end{align*}
Now assume that $\xi'\colon\EM(T')$ also satisfies property $\textbf{P}$. Then
we have 1-cells and invertible 2-cells
\begin{align*}
\varphi\colon&T\to T'\,,\qquad\phi\colon\EM(\varphi)\ohArr\xi\iso\xi'\\
\varphi'\colon&T'\to T\,,\qquad\phi'\colon\EM(\varphi')\ohArr\xi'\iso\xi
\end{align*}
and it is easy to see that there are invertible 2-cells
$\varphi\ohBmd\varphi'\iso T'$ and $\varphi'\ohBmd\varphi\iso T$, i.e. $T$
and $T'$ are equivalent. This result, however weak, is in complete agreement
with the fact that $\QQ$, as a left pseudo adjoint of $\EM$, is determined only
up to pseudo natural isomorphisms. 

On the other hand, the way we defined $\QQ$ allowed only the freedom to choose
different adjunction datas for the functors $U$, which amounts to $\QQ$
beeing unique up to 2-natural isomorphisms. Also the universality formulated
in Theorem \ref{thm: 1-univ} suggests that we should find a 2-categorical
2-adjunction generalizing it.  

Notice that the image of $\EM$ lies in a special sub-2-category of
$\LMF$ in which the 1-cells contain identity natural isomorphisms
$\kappa$. Let us call a 1-cell $\bra F,\kappa,G\ket$ strict if
$\kappa=1_{GU}=1_{VF}$. The 2-category of all objects of $\LMF$ with only
strict 1-cells between them and with all 2-cells between strict 1-cells will
be denoted by $\StMF$. 

Remember that the counit $\zeta\colon\QQ\twocom\EM\to\Bmd$ is a 2-natural
transformation. The unit $\xi\colon\LMF\to\EM\twocom\QQ$ is only pseudo
natural but the 2-cell $\Xi_\kappa$ is such that it is the identity for strict
1-cells. Therefore the restriction of $\xi$ to $\StMF$ is also 2-natural.
Denoting by $\EM_{st}$ and $\QQ^{st}$ the corresponding restricted 2-functors
we obtain an ordinary 2-adjunction
\beq
\QQ^{st}\dashv\EM_{st}\,,
\eeq
i.e., one in which the unit and counit are 2-natural transformations. Such
left adjoints $\QQ^{st}$ are already unique up to 2-natural isomorphisms. This
is reflected by the following property of the monad $\QQ(U)$ of a left
adjointable monoidal functor.
\begin{thm} \label{thm: 2-univ}
For each object $U$ of $\LMF$ there exists a bimonad $T$ and a strict
1-cell $\xi\colon U\to\EM(T)$ with the following property:
\begin{quote}
\textbf{U:} If $S$ is a bimonad and $\iota\colon U\to\EM(S)$ is a strict 1-cell
then there exists a unique monad morphism $\varphi\colon T\to S$ such that
\[
\EM(\varphi)\ohArr\xi\ =\ \iota\,.
\]
\end{quote}
If another bimonad $T'$ and another strict 1-cell $\xi'\colon U\to\EM(T')$ has
property \textbf{U} then there is a bimonad isomorphism $\psi\colon T\iso T'$
such that $\EM(\psi)\ohArr\xi=\xi'$.
\end{thm}

\section{Bialgebroids} \label{sec: bgd}

\subsection{From bialgebroids to bimonads}

Let $k$ be a commutative ring, $R$ a (possibly non-commutative) $k$-algebra. A
Takeuchi $\times_R$ bialgebra or a
left bialgebroid over $R$ in the sense of \cite{KSz} consists of 
\begin{itemize}
\item a $k$-algebra $A$
with a $k$-algebra map $s\otimes_k t\colon R\otimes_k R\op\to A$ making $A$
into an $R$-$R$ bimodule via $r\cdot a\cdot r':=s(r)t(r')a$ and 
\item a comonoid
structure $\bra A,\cop,\eps\ket$ on $A$ in $_R\M_R$  
\end{itemize}
such that
\begin{description}
\item[BGD 1.a] the image of the comultiplication $\cop(A)\subset A\otimes_R
A$ belongs to the subbimodule  
\[
A\times_R A=\{X\in A\otimes_R A\,|\,X(1\otimes
s(r))=X(t(r)\otimes 1)X,\forall r\in R\,\}
\]
which has the obvious algebra structure therefore it is meaningful to require
that 
\item[BGD 1.b] $\cop\colon A\to A\times_R A$ be a $k$-algebra map, moreover

\item[BGD 2.a] the counit $\eps$ preserves the unit,
$\eps(1_A)=1_R$  
\item[BGD 2.b] and satisfies
\[\eps(at(\eps(b)))=\eps(ab)=\eps(as(\eps(b)))\]
for all $a,b\in A$.
\end{description}

Right bialgebroids are defined analogously but using right multiplications
with $s(r)$, $t(r)$ in the definition of the $R$-$R$-bimodule structure of
$A$, so the meaning of $A\times_R A$ also changes. What is important that in
order for the category $_A\M$ of left $A$ modules to have a monoidal structure
one needs a left bialgebroid $A$ while a right bialgebroid makes $\M_A$ to be
monoidal.

Every left $A$ module $_AV$ inherits an $R$-$R$ bimodule structure via the
algebra map $s\otimes_k t$, i.e., if we denote the action of $a\in A$ on
an element $v\in V$ by $a\lact v$ then $r\cdot v\cdot
r':=s(r)t(r')\lact v$,  $r,r'\in R$. This defines the forgetful
functor $U\colon\,_A\M\to\,_R\M_R$. 

The comultiplication $\cop\colon a\mapsto a\1\o a\2$ allows to
define a monoidal product on $_A\M$ such that $U$ becomes strictly monoidal.
The monoidal product $X\oalg Y$ of the $A$-modules $X$ and $Y$ is the
$R$-$R$-bimodule $X\o_R Y$ equipped with $A$-action $a\lact(x\o y)=(a\1\lact
x)\o(a\2\lact y)$ which is well defined due to axiom (BGD 1.a) above. 

In the sequel we shall identify $R$-$R$-bimodules $X$ with left $R^e$-modules
via $(r\o r')\cdot a:=r\cdot a\cdot r'$, where $R^e:=R\o R\op$. The left
regular $A$-module $A=\,_AA$ is not only a left $R^e$-module but a right
$R^e$-module, as well. This allows to define a functor
$T:=A\o_{R^e}\under\colon \,_R\M_R\to\,_R\M_R$. 

\begin{thm} \label{thm: bmd of bgd}
Let $A$ be a left bialgebroid over $R$. 
Then the endofunctor $T=A\o_{R^e}\under$ defines a bimonad on $_R\M_R$
with structure maps
\begin{align}
\mu_X\colon &A\o_{R^e}(A\o_{R^e}X)\to A\o_{R^e} X\,,\quad a\o(b\o x)\mapsto
ab\o x \label{eq: TofA 1}\\
\eta_X\colon &X\to A\o_{R^e}X\,,\qquad x\mapsto 1_A\o x \label{eq: TofA 2}\\
\gamma_{X,Y}\colon&A\o_{R^e}(X\o_R Y)\to
                   (A\o_{R^e}X)\o_R(A\o_{R^e}Y)\,,\notag\\ 
&a\o(x\o y)\mapsto (a\1\o x)\o (a\2\o y) \label{eq: TofA 3}\\ 
\pi\colon& A\o_{R^e}R\to R\,,\qquad a\o r\mapsto
\eps(as(r))\,. \label{eq: TofA 4} 
\end{align}
\end{thm}
\begin{proof}
Since the bimodule structure of $A$ comes from $R^e$ being a subalgebra in $A$,
the monoid structure $A\o A\to A$ in $_k\M$ determines, via
the coequalizer $A\o_k A\rightarrowtail A\o_{R^e}A$, a monoid structure
$A\o_{R^e}A\to A$ in $_{R^e}\M_{R^e}$. This latter monoid structure makes $T$
into a monad on $_{R^e}\M$ with structure maps given in elementwise notation
in (\ref{eq: TofA 1}-\ref{eq: TofA 2}). Thus $\mu$ and $\eta$ satisfies the
bimonad axioms (BMD 5) and (BMD 6), i.e., $\bra T,\mu,\eta\ket$ is a monad.

The comultiplication $\cop\colon A\to A\x_R A$ defines the
comonoidality natural transformation $\gamma$ by formula (\ref{eq: TofA 3}). It
is well defined due to axiom (BGD 1.a) and it satisfies the hexagon of (BMD 1) due
to coassociativity of $\cop$. 
The other component $\pi$ of the comonoidality
structure of $T$ given in (\ref{eq: TofA 4}) is well defined due to that axiom
(BGD 2.b) implies $\eps(as(r))=\eps(at(r))$ for $a\in A$, $r\in R$. It is a
counit for $\gamma$ in the sense of the bimonad axiom (BMD 2) because $\eps$ is
the counit for $\cop$. Thus $\bra T,\gamma,\pi\ket$ is a comonoidal functor.

The compatibility condition of $\mu$ with $\gamma$ follows from the
bialgebroid axiom (BGD 1.b) while the compatibility of $\mu$ with $\pi$ follows
from the counit axiom (BGD 2.b). This proves (BMD 3). Compatibility of $\eta$ with
$\gamma$ is unitality of $\cop$ hence follows from (BGD 1.b) while compatibility
of $\eta$ with $\pi$ is the section property $\eps\c s=\id_R$. This proves
(BMD 4).
\end{proof}

\begin{rmk}
In this Section we speak about bialgebroids in the category $_k\M$ of
$k$-modules including as special cases the category $\Ab$ of Abelian groups or
the category $\Vec K$ of vector spaces over a field $K$. However, the category
$_k\M$ can be replaced with any symmetric monoidal closed category
$\bra\M,\o,i\ket$ which has coequalizers. The symmetric monoidal structure is
required to be able to speak about the monoids $R$, $R\op$ and $R^e$ while the
coequalizers are needed to define tensor products over such monoids. The
tensor product $\o_R$ becomes a monoidal product on $_R\M_R$ if $\o$ preserves
coequalizers. This latter property is guarantied if $\M$ is closed. The
definition of bialgebroids as well as the above theorem holds also in this
more general setting and should cover non-additive examples. 
\end{rmk}

\begin{exa} \label{exa: E}
The trivial left bialgebroid over $R$ is the bialgebroid $E=R\o_k R\op$
with comultiplication and counit given respectively by
\begin{align}
\cop_E&\colon E\to E\o_R E\,,\quad r\o_k r'\mapsto(r\o_k 1_R)\o_R(1_R\o_k
r')\\
\eps_E&\colon E\to R\,,\qquad r\o_k r'\mapsto rr'
\end{align}
and with source and target maps $s_E(r)=r\o_k 1_R$, $t_E(r)=1_R\o_k r$.
\end{exa}

If $T$ is the bimonad on $_R\M_R$ associated to a bialgebroid $A$ over $R$
then the category of $T$-algebras $_R\M_R^\ttT$ is monoidally isomorphic to
the category of left $A$-modules $_A\M$.

\subsection{Characterizing bimonads of bialgebroids}

The bimonad constructed in Theorem \ref{thm: bmd of bgd} is special among the
bimonads in that $T$ is $k$-linear and has a right adjoint. As a matter of
fact, the functor $\Hom(\,_{R^e}A,\under)$ maps an $R$-$R$-bimodule $X$ into
the $k$-module $\Hom(A,X)$ of $k$-linear maps $f\colon A\to X$ satisfying
$f(s(r)t(r')a)=r\cdot f(a)\cdot r'$ and equipped with $R$-$R$-bimodule
structure $r\cdot f\cdot r'=f(\under s(r)t(r'))$. This functor is the right
adjoint of $T=A\o_{R^e}\under$ with counit and unit
\begin{align}
&A\o_{R^e}\Hom(A,X)\to X\,,\qquad a\o f\mapsto f(a)\\
&X\to\Hom(A,A\o_{R^e}X)\,,\qquad x\mapsto \{a\mapsto a\o x\}\,.
\end{align}
From now on we never mention $k$-linearity although every functor on 
$k$-linear categories will be assumed $k$-linear. This means for example that
bimonads on $k$-linear categories will be assumed to have $k$-linear
underlying functors. Let $\Bmd_k$ be the 2-category of such bimonads.

In this Subsection we will show that the single property of
having a right adjoint already characterizes the bialgebroids within the
bimonads of $\Bmd_k$. The summary is this.
\begin{thm} \label{thm: bgd bimonads} 
Let $R$ be a monoid in $_k\M$ and $\bra
T,\gamma,\pi,\mu,\eta\ket$ be a bimonad on $_R\M_R$. Then $T$ is isomorphic
to the bimonad of a bialgebroid over $R$ if and only if it has a right
adjoint.  
\end{thm} 
Only sufficiency requires a proof. Nevertheless we will
give a detailed proof divided into a series of Lemmas that contain both
necessary and sufficient conditions. 
\begin{lem} \label{lem: 1} 
Let $E$ be a $k$-algebra and $T\colon
\,_E\M\to\,_E\M$ be an endofunctor on the category of left $E$-modules. Then
there exists an $E$-$E$-bimodule $M$ such that $T\cong M\o_E\under$ if and
only if $T$ has a right adjoint. 
\end{lem}
\begin{proof}
Necessity: $M\o_E\under$ has a right adjoint
$\hom(M,\under):=\Hom_E(\,_EM,\under)$ inheriting its left $E$ module
structure from the right $E$-action on $M$. The adjunction relation
\begin{align}
\Hom_E(X,\hom(M,Y))&\iso\Hom_E(M\o_E X,Y)\\
\notag f&\mapsto \{m\o x\mapsto f(x)(m)\}
\end{align}
for left $E$-modules $X$ and $Y$ is a standard hom-tensor relation.

Sufficiency: Let $H$ be a right adjoint to $T$. Then, $M:=T(\,_EE)$ being an
$E$-$E$-bimodule via $E\op\cong\End(\,_EE)$ as well as $E$, we have 
\begin{align*}
\Hom_E(TX,Y)&\cong\Hom_E(X,HY)\cong\Hom_E(E\o_E X,HY)\\
&\cong\Hom_E(X,\hom(E,HY))\\
&\cong\Hom_E(X,\hom(TE,Y))\\
&\cong\Hom_E(TE\o_E X, Y)
\end{align*}
implying $TX\cong TE\o_E X$.
\end{proof}

\begin{lem} \label{lem: 2}
Let $E$ be as before and $\bra T,\mu,\eta\ket$ be a monad on $_E\M$. Then
there exists a monoid $A$ in $_E\M_E$ and a monad
isomorphism $A\o_E\under\cong T$ if and only if $T$ has a right adjoint. 
\end{lem}
\begin{proof}
Necessity: This is the same as the necessity part of the previous Lemma.

Sufficiency: By the previous Lemma there is a bimodule $A$ and an isomorphism
$\nu\colon T\iso T':=A\o_E\under$ of functors. The natural transformations
$\mu$ and $\eta$ can be passed to $T'$ via $\nu$ to get a monoid $\bra
T',\mu',\nu'\ket$. Since the powers of $T'$ have hom-functors as right
adjoints and the natural transformations between them -- by the Yoneda Lemma
-- arise from bimodule maps between the tensor powers of $A$, it is easy to see
that the monad structure on $T'$ is that of arising from a monoid structure on
$A$. 
\end{proof}

The next Lemma provides an important class of examples of lax monoidal
functors.
\begin{lem} \label{lem: 3}
Let $\bra\M,\o,i\ket$ be a monoidal category with coequalizers and assume that
$x\o\under$ and $\under\o x$ preserve coequalizers for all objects $x$ of $\M$.
Then for any monoid $\bra R,\nu,\iota\ket$ in $\M$ the $R$-$R$-bimodules
$X=\bra x,\lambda_X\colon R\o x\to x,\rho_X\colon x\o R\to x\ket$ in $\M$
form a monoidal category $_R\M_R$ with monoidal product $\o_R$ arising from a
choice of coequalizers $x\o(R\o y)\rightrightarrows x\o y\rightarrowtail x\o_R
y$ for each pair of objects. The forgetful functor $\Phi\colon\,_R\M_R\to\M$
mapping $\bra x,\lambda_X,\rho_X\ket$ to $x$ is lax monoidal with
\beq
\Phi_{X,Y}\colon \Phi X\o\Phi Y\to\Phi(X\o_R Y)
\eeq
being the chosen coequalizer and 
\beq
\Phi_0\colon i\to\Phi R
\eeq
being just the unit $\iota$ of the monoid $R$.
\end{lem}
\begin{proof}
The proof is standard and therefore omitted. In case of $\M=\,_k\M$
(for which the notation $_k\Ab$ would be more logical) and $R$ a $k$-algebra
the statement is definitely common lore.
\end{proof}

Combining the results of the last two Lemma with the fact that a lax monoidal
functor $\Phi\colon \B\to\M$ maps monoids $\bra A,\mu,\eta\ket$ in $\B$ into
monoids 
\[
\bra \Phi A,\Phi(\mu)\c\Phi_{A,A},\Phi(\eta)\c\Phi_0\ket
\]
in $\M$, we immediately obtain
\begin{cor} \label{cor: 4}
Let $E$ be a $k$-algebra and $\bra T,\mu,\eta\ket$ be a monad on $_E\M$.
Then there exists a $k$-algebra extension $E\to A$ and a monad isomorphism
$A\o_E\under\iso T$ if and only if $T$ has a right adjoint.
\end{cor}

Now we investigate the coalgebra properties of the $k$-algebra $A$.
\begin{lem} \label{lem: 5}
Let $R$ be a $k$-algebra, $E=R\o_k R\op$ be its enveloping algebra and let
$\bra T,T^2,T^0\ket$ be a lax comonoidal endofunctor on the monoidal category
$_E\M$. Then there exists 
\begin{itemize}
\item an $E$-$E$-bimodule $A$,
\item a comonoid $\bra A,\cop,\eps\ket$ in $_E\M$ in which the
$R$-$R$-bimodule structure of $A$ comes from $_EA$,  
\item and a comonoidal natural isomorphism $A\o_E\under\iso T$
\end{itemize}
if and only if $T$ has a right adjoint. 
\end{lem} 
\begin{proof}
Necessity: This holds by Lemma \ref{lem: 1} even without comonoid structure.

Sufficiency: If $T$ has a right adjoint then Lemma \ref{lem: 1} ensures the
existence of a bimodule $A$ in $_E\M_E$ and a natural isomorphism
$A\o_E\under\iso T$. Using this isomorphism we can put a comonoidal structure
on $A\o_E\under$ making the isomorphism into a comonoidal natural isomorphism.
Let 
\begin{align*}
\gamma_{X,Y}&\colon A\o_E(X\o_R Y)\to (A\o_E X)\o_R(A\o_E Y)\\
\pi&\colon A\o_E R\to R
\end{align*}
be the lax comonoidal structure we obtained that way. 
Since $E$ is a generator for $_E\M$ and $\gamma$ is natural, the
components $\gamma_{X,Y}$ are completely determined by $\gamma_{E,E}$. As a
matter of fact, for $x\in X$ let $f_x\colon E\to X$ be defined by $f(r\o
r'):=r\cdot x\cdot r'$. Similarly, let $g_y$ be the analogue for $Y$. Then for
all $x\in X$, $y\in Y$ and $a\in A$ 
\[
\gamma_{X,Y}(a\o_E(x\o_R y))=[(A\o_E f_x)\o_R(A\o_E g_y)]\c
\gamma_{E,E}(a\o_E(1_E\o_R 1_E))
\]
or, introducing
\begin{align}
\cop(a)&:=(\rho_A^{-1}\o_R\rho_A^{-1})\c\gamma_{E,E}(a\o_E(1_E\o_R 1_E))\\
&=a\1\o_R a\2\ \in A\o_R A
\end{align}
where $\rho_A\colon A\o_E E\iso A$ denotes the obvious isomorphism, we obtain
formula (\ref{eq: TofA 3}). 
Inserting this expression of $\gamma$ into the diagrams (\ref{dia: BMD 2}) 
and take the special case $x=E$ one obtains 
\[
\pi(a\1\o_E 1_R)\cdot a\2\ =\ a\ =\ a\1\cdot\pi(a\2\o_E 1_R)\,,\qquad a\in
A\,.
\]
Therefore $\cop$ is counital with counit 
\beq
\eps(a):=\pi(a\o_E 1_R)
\eeq
and it is left for an exercise to show that (\ref{dia: BMD 1}) implies that
$\cop$ is coassociative.
\end{proof}
\begin{rmk} \label{rmk: cop eps}
The comultiplication and counit of $A$ can be expressed in terms of the
coring structure of $E$ of Example \ref{exa: E} and in terms of the comonoidal
structure (\ref{eq: TofA 3}-\ref{eq: TofA 4}) of $T=A\o_E\under$ as follows. 
The comultiplication $\cop$ is the composite 
\beq
\begin{CD}
A@>\rho_A^{-1}>\sim>A\o_E E@>T\cop_E>>A\o_E(E\o_R E)@>\gamma_{E,E}>>(A\o_E
E)\o_R(A\o_E E)\\
@.@.@.@V{\wr}V{\rho_A\o_R\rho_A}V\\
@.@.@.A\o_R A
\end{CD}
\eeq
while the counit is the composite
\beq
\begin{CD}
A@>\rho_A^{-1}>\sim>A\o_E E@>T\eps_E>>A\o_E R@>\pi>>R
\end{CD}
\eeq
in the category of $R$-$R$-bimodules.
\end{rmk}

It is interesting that the Takesaki $\x_R$ product appears naturally already
in the bimodule context, i.e., without the algebra structures, as the next
Lemma shows.
\begin{lem} \label{lem: 6}
Let $R$ and $E$ be as in Lemma \ref{lem: 5} and let $\bra A,\cop,\eps\ket$ be a
comonoid in $_R\M_R\equiv\,_E\M$ for some $E$-$E$ bimodule $A$ such that the
endofunctor $A\o_E\under$ on $_E\M$ is lax comonoidal. Then $\cop$ and $\eps$
satisfy
\begin{align}
\a\1\cdot t_E(r)\o_R a\2&=a\1\o_R a\2\cdot s_E(r) \label{eq: Tak cop}\\
\eps(a\cdot t_E(r))&=\eps(a\cdot s_E(r)) \label{eq: Tak eps}
\end{align}
for all $a\in A$, $r\in R$.
\end{lem}
\begin{proof}
The proof uses essentially that $\cop$ and $\eps$ can be expressed in terms of
$\gamma_{E,E}$ and $\pi$, see Remark \ref{rmk: cop eps}. First of all, the
right action $\under\cdot e:=\rho_A(\under\o_E e)$ of an element of $E$ on $A$
commutes with the left $E$ action therefore it is an $R$-$R$ bimodule map.
For a fixed $a\in A$ choose a finite set of $a_{ij},b_{ik}\in A$ and
$e_j,f_k\in E$ such that 
\[
\gamma_{E,E}(a\o_E(1_E\o_R 1_E))=(a_{ij}\o_E e_j)\o_R(b_{ik}\o_E f_k)
\]
with summations understood. Then applying naturality of $\gamma_{E,E}$ twice
for any $r\in R$ we can write
\begin{align*}
&(a_{ij}\o_E e_j t_E(r))\o_R(b_{ik}\o_E f_k)=\gamma_{E,E}(a\o_E(t_E(r)\o_R
1_E))\\
&\gamma_{E,E}(a\o_E(1_E\o_R s_E(r)))=(a_{ij}\o_E e_j)\o_R(b_{ik}\o_E f_k
s_E(r)) 
\end{align*}
implying that $\cop(a)=a_{ij}\cdot e_j\o_R b_{ik}\cdot f_k$ satisfies
(\ref{eq: Tak cop}).
In order to get (\ref{eq: Tak eps}) use  
\begin{align*}
\pi(a\cdot t_E(r)\o_E 1_R)&=\pi(a\o_E t_E(r)\cdot 1_r)=\pi(a\o_E r)
=\pi(a\o_E s_E(r)\cdot 1_R)\\
&=\pi(a\cdot s_E(r)\o_E 1_R)\,.
\end{align*}
\end{proof}

Now we can finish the proof of Theorem \ref{thm: bgd bimonads} as follows.
\begin{proof}
That the bimonad of a bialgebroid has a right adjoint was shown at the
beginning of this subsection. Assume $T$ is a bimonad on $_R\M_R$ with a right
adjoint. Then by Corollary \ref{cor: 4} there is an algebra extension $A$ of
$R^e$ and a monad isomorphism $\nu\colon T\iso A\o_E\under$. Use this $\nu$
to pass the bimonad structure of $T$ to the functor $A\o_E\under$. Then $\nu$
becomes a bimonad isomorphism. Now $A\o_E\under$ has a right adjoint
therefore by Lemma \ref{lem: 5} there is a comonoidal natural
isomorphism $\chi\colon A\o_E\under\iso B\o_E\under$ for some $R$-coring
and $E$-$E$ bimodule $B$. This latter isomorphism determines an $E$-$E$
bimodule isomorphism $A\iso B$ which can be used to make $A$ into a comonoid
in $_E\M$. Now the bimonad $A\o_E\under$ has structure maps as in (\ref{eq:
TofA 1}-\ref{eq: TofA 4}) in which $\cop$ and $\eps$ give rise to an
$R$-coring structure on $A$ and satisfy 
\begin{align}
\cop(a)(t(r)\o 1_A)&=\cop(a)(1_A\o_R s(r))\\
\eps(as(r))&=\eps(at(r))
\end{align}
by Lemma \ref{lem: 6}. It remains to use the bimonad axioms (BMD 3) and (BMD
4). Inserting (\ref{eq: TofA 3}) into the first diagrams of (\ref{dia: BMD 3})
and (\ref{dia: BMD 4}),
after a little calculation one obtains that $\cop$ has to be multiplicative
and unit preserving, respectively. Substituting (\ref{eq: TofA 4}) into the
remaining diagrams of (\ref{dia: BMD 3}) and (\ref{dia: BMD 4}) one immediately
arrives to the two bialgebroid axiom (BGD 2.a) and (BGD 2.b).
This proves that $A$ is a left bialgebroid over $R$ and that its bimonad is
isomorphic to $T$ via a comonoidal natural isomorphism.
\end{proof}

\subsection{Tannaka reconstruction for bialgebroids}

By characterizing bialgebroids as the bimonads with right adjoints (on bimodule
categories) a natural definition arises for what to be the bialgebroid
morphisms and transformations.
\begin{defi} \label{def: bgd}
Let $\Bgd_k$ denote the 1-full and 2-full sub-2-category of $\Bmd_k$ the
objects of which are the $k$-linear bimonads $T\colon \T\to\T$ with
right adjoint where $\T$ is isomorphic to $_R\M_R$ for some $k$-algebra $R$.
The objects of the form $A\o_{R^e}\under$ for some bialgebroid $A$ over $R$
are called proper bialgebroids. 
\end{defi}

Note that we could define $\Bgd_k$ to be 2-replete and not only 1-replete by
allowing for objects all bimonads that are equivalent to proper bialgebroid
bimonads. Still the above definition works well with the Eilenberg-Moore
construction.
\begin{pro}
Let $T\colon \T\to\T$ be an object in $\Bgd_k$. Then
$\EM(T)=U^\ttT\colon$ $\T^\ttT\to \T$ is a monoidal functor with both
left and right adjoints.
\end{pro}
\begin{proof}
Strict monoidality of $U^\ttT$ and existence of left adjoint follows from
Propositions \ref{pro: mon EM} and \ref{pro: F^T}. Existence of right adjoint
follows from \cite[Corollary V.8.3]{MacLane-Moerdijk}.
\end{proof}

The converse of the above proposition, namely that a monoidal functor
$U\colon\C\to\T$ with both left and right adjoints determines a
bimonad $\ttT=\QQ(U)$ with right adjoint is obvious since the underlying
functor is now a product of two functors $T=UF$ with both $U$ and $F$ having a
right adjoint. Therefore appropriate restrictions of the 2-functors $\EM$ and
$\QQ$, denoted by the same letters, provide an adjunction and a Tannakian
theory for bialgebroids.

In the following definition $\LMF_k$ denotes the $k$-linear version of $\LMF$
of Section \ref{sec: 2func} with only strict 1-cells.
\begin{defi}
Let $\AMF_k$ be the 1-full and 2-full sub-2-category of $\LMF_k$ the objects of
which are the monoidal functors $U\colon\C\to\T$ with both left and right
adjoints and with target category $\T$ that is isomorphic to some bimodule
category over $_k\M$.  
\end{defi}
It follows that the restrictions of $\QQ$ and $\EM$ to 2-functors
$\AMF_k\to\Bgd_k$ and $\Bgd_k\to\AMF_k$, respectively, constitute a
2-adjunction $\QQ\dashv\EM$.

The following corollary is a direct consequence of Theorem \ref{thm: 2-univ}
and of the above definitions.
\begin{cor}
Let $R$ be a $k$-algebra and $U\colon \C\to\,_R\M_R$ be a $k$-linear monoidal
functor with left and right adjoints. Then there is a bialgebroid $A$ over $R$
and a monoidal functor $K\colon \C\to\,_A\M$ such that 
\begin{enumerate}
\item $U$ factorizes as $U=U_A\, K$ through the strict monoidal forgetful
functor $U_A$ of the bialgebroid,
\item if $B$ is another bialgebroid over some $k$-algebra $S$ such that
there exist monoidal functors $F\colon \C\to\,_B\M$ and $G\colon
\,_R\M_R\to\,_S\M_S$ satisfying $GU=U_B\, F$ then there exists a unique
bialgebroid morphism $\bra G,\varphi\ket\colon A\to B$, i.e., a unique
ambimonoidal natural transformation
\[
\varphi\colon B\o_{S^e}G(\under)\to G(A\o_{R^e}\under)
\]
such that $G^\varphi K=F$.
\end{enumerate}
The bialgebroid $A$ with the above properties is unique up to isomorphisms.
\end{cor}

The following representation theorem for bialgebroids, in turn, is a
consequence of Theorem \ref{thm: rep}.
\begin{cor}
Let $\C$ be a $k$-linear monoidal category and $R$ be a $k$ algebra.
Then for a $k$-linear functor $U\colon\C\to\,_R\M_R$ the following conditions
are equivalent: 
\begin{enumerate}
\item There exists a bialgebroid $A$ over $R$ and a $k$-linear monoidal
category equivalence $K\colon\C\to \,_A\M$ such that $U_A\, K=U$.
\item $U$ is monadic, monoidal and has a right adjoint.
\end{enumerate}
\end{cor}

\subsection{Bialgebroid maps}

According to Definition \ref{def: bgd} the morphisms $\bra G,\varphi\ket$
from a bialgebroid $A$ over $R$ to another $B$ over $S$ consists of a lax
monoidal functor $G\colon\,_R\M_R\to\,_S\M_S$ and of an ambimonoidal natural
transformation $\varphi\colon B\o_{S^e}G(\under)\to G(A\o_{R^e}\under)$
satisfying the two diagrams (\ref{dia: mnd mor}). These conditions are rather
complicated for a general functor $G$ so we can only give some special
examples. The simplest are the bialgebroid maps.

Assume that $G$ arises from a $k$ algebra homomorphism $\omega\colon S\to R$,
i.e., $G$ is a lax monoidal forgetful functor 
\beq
G=\Phi^\omega\colon\ _R\M_R\to\ _S\M_S\,,\qquad _RX_R
\mapsto\ _{\omega(S)}X_{\omega(S)}
\eeq
similar to the $\Phi$ of Lemma \ref{lem: 3}. In this case the natural
transformation 
\[
\varphi_X\colon B\o_{S^e}\Phi^\omega(X)\to \Phi^\omega(A\o_{R^e}X) 
\]
is completely determined by 
\[
\varphi_{R^e}\colon B\o_{S^e}\,X\to A\o_{R^e}X
\]
since $R^e$ is a generator. Naturality of $\varphi_{R^e}$ alone in turn gives
\[
\varphi_{R^e}(b\o_{S^e}(r\o r'))=\varphi(b)\o_{R^e}(r\o r')
\]
where the $S$-$S$-bimodule map $\varphi\colon B\to \Phi^\omega(A)$ is the
composite 
\beq \label{dia: fi}
\begin{CD}
B\iso B\o_{S^e} S^e@>B\o_{S^e}(\omega\o \omega')>>B\o_{S^e}R^e
@>\varphi_{R^e}>>\Phi^\omega(A)\o_{R^e}R^e\iso\Phi^\omega(A)
\end{CD}
\eeq
where note that in $\Phi^\omega(A)$ only the left $R^e$ action is forgotten,
the right one is intact. Inserting the expression $\varphi_X(b\o
x)=\varphi(b)\o x$ into the monad morphism axioms (\ref{dia: mnd mor}) we
obtain that $\varphi\colon B\to A$ is an algebra map.
Since it is also an $S$-$S$ bimodule map by its definition (\ref{dia: fi}),
we obtain the identities
\begin{align}
\varphi\c s_B&=s_A\c\omega \label{eq: bgd map 1}\\
\varphi\c t_B&=t_A\c\omega \label{eq: bgd map 2}
\end{align}
Now inserting to the
ambimonoidality axioms of Definition \ref{def: ambi} we obtain that
$\varphi\colon B\to A$ preserves the coalgebra structure in the sense of the
equations 
\begin{align}
\Phi^\omega(\cop_A)\c\varphi&=\Phi^\omega_{A,A}\c(\varphi\o_S\varphi)\c\cop_B
\label{eq: bgd map 3}\\
\Phi^\omega(\eps_A)\c\varphi&=\Phi^\omega_0\c\eps_B
\label{eq: bgd map 4}
\end{align}
or in Hopf algebraist notation
\begin{align}
\varphi(b)\1\o_R\varphi(b)\2&=\varphi(b\1)\o_R\varphi(b\2)\\
\eps_A(\varphi(b))&=\omega(\eps_B(b))
\end{align}
for all $b\in B$. The equations (\ref{eq: bgd map 1}-\ref{eq: bgd map 2}-
\ref{eq: bgd map 3}-\ref{eq: bgd map 4}) define what is called a bialgebroid
map in \cite{Strbg}. A bialgebroid map is completely determined by
$\varphi\colon B\to A$ since $\omega=\eps_A\c\varphi\c s_B$.

\subsection{Bimodule induced bialgebroid morphisms}

Another class of bialgebroid morphisms are obtained if we take the
functor $G\colon\,_R\M_R\to\,_S\M_S$ to be $GX=G\o_{R^e}X$ for some
$S^e$-$R^e$-bimodule $G$. The natural transformation $\varphi\colon 
B\o_{S^e}G(\under)\to G(A\o_{R^e}\under)$ then becomes expressed in terms of a
bimodule map 
\beq
\varphi\colon B\o_{S^e}G\to G\o_{R^e}A\ \in\ _{S^e}\M_{R^e}
\eeq
as $\varphi_X(b\o g\o x)=\varphi(b\o g)\o x$. If we insert this expression
into the two monad morphism diagram (\ref{dia: mnd mor}) and into the two
ambimonoidality diagram (\ref{dia: ambi1}-\ref{dia: ambi2}) we obtain four
relations between $G$ and $\varphi$ that are reminiscent of the entwining
structure of Brzezi\'nski and Majid \cite{Brz-Maj}, although not the same. 

At first notice that lax monoidality of the functor $\bra G, G_2,G_0\ket$
imposes a monoid structure on the bimodule $G$ in the category of
$S$-$S$-bimodules but also satisfies dual analogues of the bialgebroid
comultiplication property (BGD 1.b) from the right hand side due to naturality
of $G_2$. The monad morphism axioms imply the following two conditions
\begin{align}
\varphi\c(\mu^B\o_{S^e}G)&=(G\o_{R^e}\mu^A)\c(\varphi\o_{R^e}A)\c
(B\o_{S^e}\varphi)\\
\varphi\c(\eta^B\o_{S^e}G)&=G\o_{R^e}\eta^A
\end{align}
where we use $\bra A,\mu^A,\eta^A\ket$ to denote the algebra structure of $A$
in $_{R^e}\M_{R^e}$ and similarly for $B$. The ambimonoidality conditions,
however, are not so easy to formulate only in terms of the coproducts
$\cop_A$ and $\cop_B$ and not the natural transformations they define via
(\ref{eq: TofA 3}). So let us specialize ourselves to $R=S=k$ and and
assume that on each bimodule $G$, $A$, and $B$ all the $k$-actions coincide.
So we have $k$-bialgebras $A$ and $B$, a $k$-algebra $G$ and a $k$-linear map
$\varphi\colon B\o G\to G\o A$. Then we can use the notation
\beq
\varphi(b\o g)\ =\ g_\varphi\o b^\varphi\ \in G\o A
\eeq
and write all the four conditions in a simple way
\begin{align}
g_\varphi\o (bb')^\varphi&=g_{\varphi'\varphi}\o b^\varphi b^{\varphi'}\\
g_\varphi\o (1_B)^\varphi&=g\o 1_A\\
g_\varphi g'_{\varphi'}\o b^\varphi\1\o b^{\varphi'}\2&=
(gg')_\varphi\o(b^\varphi)\1\o(b^\varphi)\2\\
(1_G)_\varphi\eps_A(b^\varphi)&=1_G\eps_B(b)
\end{align}

Returning to the case of general $R$ and $S$ the bialgebroid morphism $\bra
G\o_{R^e}\under,\varphi\ket$ described above is in fact the most general
possible if we require it to be an equivalence of the objects $A$ and $B$ in
$k$-$\Bgd$. This follows from Morita theory since $_R\M_R\rarr{G}\,_S\M_S$
should be an equivalence and using the isomorphisms $_{R^e}\M\equiv\,_R\M_R$
and $_{S^e}\M\equiv\,_S\M_S$ the $G$ has to be naturally isomorphic to a
functor $G\o_{R^e}\under$ with a Morita equivalence bimodule 
$_{S^e}G_{R^e}$. That is to say, the rings $R$ and $S$ are
$\sqrt{\text{Morita}}$-equivalent \cite{Takeuchi: sqrt M}. Ordinary Morita
equivalence $R\sim S$ arises under the further assumption that $_{S^e}G_{R^e}$
is the $k$-tensor  product of equivalence bimodules $_RH_S$ and $_SH'_R$. This
latter situation is the Morita base change proposed by Schauenburg
\cite{Schauenburg: Morita} while the former was named as
$\sqrt{\text{Morita}}$ base change.

\subsection{An exotic example: Hom}

For the tired Reader's sake let stand here an example of a bimonad that is not
a bialgebroid. It shows that every set is a bimonad in a canonical way.

Let $\Set$ be the category of small sets equipped with the Cartesian closed
monoidal structure $\bra\Set,\x,1\ket$ with some one element set $1=\{\star\}$.
Every object $C$ in $\Set$ is a comonoid in a unique way, namely by the
diagonal mapping $\cop_C\colon x\mapsto \bra x,x\ket$ and by the constant
mapping $\eps_C\colon C\to 1$. This comonoid structure makes the endofunctor
$T:=\Hom(C,\under)$ into a monad with multiplication and unit 
\begin{align}
\mu_A&\colon T^2A\to TA\,,\qquad \mu_A(f)(c)=f(c)(c)\\
\eta_A&\colon A\to TA\,,\qquad \eta_A(a)(c)=a
\end{align}
respectively. Moreover, $T$ is comonoidal with 
\begin{align}
\gamma_{A,B}&\colon T(A\x B)\iso T(A)\x T(B)\,,\quad f\mapsto\bra p_1\c f,p_2\c
f\ket\\
\pi&\colon T(1)\iso 1\,, f\mapsto \star
\end{align}
where $p_i$ are the projections of the product $A\x B$.
Now it is an easy exercise to check that both $\mu$ and $\eta$ are comonoidal
natural transformations, so $\bra T,\gamma,\pi,\mu,\eta\ket$ is a bimonad.

A $T$-algebra for this bimonad is a set $A$ and a function
$\alpha\colon\Hom(C,A)\to A$ such that the two diagrams of (\ref{eq: T-alg})
commute. For a finite set $C$ with $n$ elements such a function $\alpha$ can be
identified with an $n$ variable function on $A$ with values in $A$. Then the
$T$-algebra conditions become the following equations for $\alpha$.
\begin{align}
&\alpha(\alpha(a_{11},\dots,a_{1n}),\dots,\alpha(a_{n1},\dots,a_{nn}))=
\alpha(a_{11},\dots,a_{nn})\\
&\alpha(a,\dots,a)=a
\end{align}
for all $a_{ij}\in A$ and $a\in A$. There are solutions that are evaluations
at an element of $C$, let's say, $\alpha(a_1,\dots,a_n)=a_i$. But there are
solutions that are not evaluations, the free $T$-algebras for example. A free
$T$-algebra $\bra \Hom(C,A),\mu_A\ket$ is a product set $\Pi^n A$ with action
\beq
\alpha(\bra a_{11},\dots,a_{1n}\ket,\dots,\bra a_{n1},\dots,a_{nn}\ket)
=\bra a_{11}\dots,a_{nn}\ket
\eeq
Of course, the solutions form a monoidal category $\Set^\ttT$ by Proposition
\ref{pro: mon EM}, otherwise the general solution for $T$-algebras is not known
to the author.

\end{document}